\documentclass[final]{siamltex}

\newtheorem{assumption}{Assumption}[section]

\renewcommand{\div}{\mbox{div}}
\newcommand{\x}{{\mathbf x}}
\newcommand{\y}{{\mathbf y}}
\newcommand{\z}{{\mathbf z}}
\newcommand{\w}{{\mathbf w}}

\newcommand{\uu}{{\mathbf u}}
\newcommand{\vv}{{\mathbf v}}

\newcommand{\n}{{\mathbf n}}
\renewcommand{\mathbf}[1]{\mbox{\boldmath$#1$}}
\newcommand{\V}{\mathbf V}
\newcommand{\W}{\mathbf W}
\newcommand{\UU}{\mathbf U}
\newcommand{\F}{\mathbf F}

\title{Homogenized models for a short-time filtration  and
 for acoustic waves propagation in a porous  media}

\author{Anvarbek Meirmanov\thanks{Math. Dept., Belgorod
State University, ul.Pobedi 85, 308015 Belgorod, Russia ({\tt
meirmanov@bsu.edu.ru}).}}

\begin{document}

\maketitle

\begin{abstract}
We consider a linear system of differential equations describing a
joint motion of elastic porous body and fluid occupying porous
space.  The rigorous justification, under various conditions imposed
on physical parameters, is fulfilled for homogenization procedures
as the dimensionless size of the pores tends to zero, while the
porous body is geometrically periodic and a characteristic time of
processes is small enough. Such kind of models may describe, for
example, hydraulic fracturing or  acoustic or seismic waves
propagation. As the results, we derive different types of
homogenized equations involving non-isotropic Stokes system for
fluid velocity coupled with acoustic equations for the solid
component or different types of acoustic equations, depending on
ratios between physical parameters. The proofs are based on
Nguetseng's two-scale convergence method of homogenization in
periodic structures.
\end{abstract}

\begin{keywords}
Stokes equations, Lam\'{e}'s equations, wave equation, hydraulic
fracturing, two-scale convergence, homogenization of periodic
structures
\end{keywords}

\begin{AMS}
 35M20, 74F10, 76S05
\end{AMS}

\pagestyle{myheadings} \thispagestyle{plain}
\markboth{A.Meirmanov}{HOMOGENIZED MODELS FOR FILTRATION AND
ACOUSTIC}

\section{Introduction}
The paper addresses the problem of a joint motion of a deformable
solid (\emph{elastic skeleton}), perforated by  system of channels
or pores (\emph{porous space}) and a fluid, occupying  porous space.
In dimensionless variables (without primes)

$$ \x'=L \x,\quad t'=\tau t,\quad \w'=L \w,
\quad \rho'_s= \rho_0 \rho_s,\quad \rho'_f =\rho_0 \rho_f,\quad
\F'=g\F,$$ differential equations of the problem in a domain $\Omega
\in \textbf{R}^{3}$ for the dimensionless displacement vector $\w$
of the continuum medium  have a form:
\begin{eqnarray} \label{0.1}
& \displaystyle \alpha_\tau \bar{\rho} \frac{\partial^2
\w}{\partial t^2}=\div_x \textbf{P} +\bar{\rho}\F,\\
 \label{0.2} & \displaystyle \textbf{P} =
\bar{\chi}\alpha_\mu D\Bigl(\x,\frac{\partial \w}{\partial t}\Bigr)
+(1-\bar{\chi})\alpha_\lambda D(x,\w)-(q+\pi)I ,\\
 \label{0.3}
& \displaystyle q=p+\frac{\alpha_\nu}{\alpha_p}\frac{\partial
p}{\partial t},\\
 \label{0.4}
& \displaystyle p+\bar{\chi} \alpha_p \div_x \w=0,\\
\label{0.5} & \displaystyle \pi+(1-\bar{\chi}) \alpha_\eta \div_x
\w=0.
\end{eqnarray}

 The problem is endowed   with homogeneous initial and boundary conditions
\begin{equation} \label{0.6}
\w|_{t=0}=0,\quad \frac{\partial \w}{\partial t}|_{t=0}=0,\quad
\x\in \Omega
\end{equation}
\begin{equation} \label{0.7}
\w=0,\quad  \x \in S=\partial \Omega, \quad t\geq 0.
\end{equation}
Here and further we use  notations
 $$ D(x,\uu)=(1/2)\left(\nabla_x \uu +(\nabla_x \uu)^T\right),
 \quad \bar{\rho}=\bar{\chi}\rho_f +(1-\bar{\chi})\rho_s.$$
 In
this model the characteristic function of the porous space
$\bar{\chi}(\x)$ and a dimensionless vector $\F(\x,t)$  of
distributed mass forces are known functions.

Dimensionless constants $\alpha_i$ $(i=\tau,\nu,\ldots)$ are defined
by the formulas $$ \alpha_\tau =\frac{ L}{g \tau^2}, \quad
\alpha_\nu = \frac{\nu }{\tau Lg\rho_0}, \quad \alpha_\mu
=\frac{2\mu }{\tau Lg\rho_0},$$
 $$ \alpha_p =\frac{c^{2}\rho _{f}}{ Lg},
 \quad \alpha_\eta =\frac{\eta}{ Lg\rho_0},
 \quad \alpha_\lambda =\frac{2 \lambda }{ Lg\rho_0},$$
where $\mu$ is the  viscosity of fluid or gas, $\nu$ is the bulk
viscosity of fluid or gas, $\lambda$ and $\eta$ are elastic
Lam\'{e}'s constants, $c$ is a  speed of sound in fluid, $L$ is a
characteristic size of the domain in consideration and $\tau$ is a
characteristic time of the process.

For more details about  Eqs.(\ref{0.1})-- (\ref{0.5})  see
\cite{AM}, \cite{B-K}, \cite{S-P}.

Our aim is to derive all possible limiting regimes (homogenized
equations) for the problem (\ref{0.1})-- (\ref{0.7}) as
$\varepsilon\searrow 0$.

 To do that we accept the following constraints
\begin{assumption} \label{assumption1}
domain  $\Omega =(0,1)^3$ is a periodic repetition of an elementary
cell  $Y^\varepsilon =\varepsilon Y$, where $Y=(0,1)^3$ and quantity
$1/\varepsilon$ is integer, so that $\Omega$ always contains an
integer number of elementary cells $Y_i^\varepsilon$. Let $Y_s$  be
a "solid part" of $Y$, and the "liquid part"  $Y_f$ -- is its open
complement. We denote as $\gamma = \partial Y_f \cap \partial Y_s$
and $\gamma $ is $C^{1}$-surface.
 A porous space  $\Omega ^{\varepsilon}_{f}$  is the periodic repetition of
the elementary cell $\varepsilon Y_f$, and solid skeleton  $\Omega
^{\varepsilon}_{s}$  is the periodic repetition of the elementary
cell $\varepsilon Y_s$. A boundary  $\Gamma^\varepsilon =\partial
\Omega_s^\varepsilon \cap \partial \Omega_f^\varepsilon$  is the
periodic repetition in  $\Omega$ of the boundary $\varepsilon
\gamma$. \textbf{The "solid skeleton"  $\Omega _{s}^\varepsilon$ and
the "porous space" $\Omega ^{\varepsilon}_{f}$ are connected
domains}.
\end{assumption}
In these assumptions

$$\bar{\chi}(\x)=\chi^{\varepsilon}(\x)=\chi
 \left(\x / \varepsilon\right),$$
$$\bar{\rho}=\rho^{\varepsilon}(\x)=\chi^{\varepsilon}(\x)\rho _{f}+
(1-\chi^{\varepsilon}(\x))\rho_{s}.$$

Let $\varepsilon$ be a characteristic size of pores $l$ divided by
the characteristic size $L$ of the entire porous body:
$$\varepsilon =\frac{l}{L}.$$

Suppose that all dimensionless parameters depend on the small
parameter $\varepsilon$ and there exist limits (finite or infinite)
$$\lim_{\varepsilon\searrow 0} \alpha_\mu(\varepsilon) =\mu_0, \quad
\lim_{\varepsilon\searrow 0} \alpha_\lambda(\varepsilon) =\lambda_0,
\quad \lim_{\varepsilon\searrow 0}
\alpha_\tau(\varepsilon)=\tau_{0},$$
$$\lim_{\varepsilon\searrow 0} \alpha_\eta(\varepsilon) =\eta_0,
\quad \lim_{\varepsilon\searrow 0} \alpha_p(\varepsilon) =p_{*},
\quad \lim_{\varepsilon\searrow 0}\alpha_\nu(\varepsilon) =\nu_0,$$
$$\lim_{\varepsilon\searrow 0}\frac{\alpha_\mu}{\varepsilon^{2}}=\mu_{1},
\quad \lim_{\varepsilon\searrow
0}\frac{\alpha_\lambda}{\varepsilon^{2}}=\lambda_{1}.$$

The first research with the aim of finding limiting regimes in the
case when the skeleton was assumed to be an absolutely rigid body
was carried out by E. Sanchez-Palencia and L. Tartar.  E.
Sanchez-Palencia \cite[Sec. 7.2]{S-P} formally obtained Darcy's law
of filtration using the method of  two-scale asymptotic expansions,
and L. Tartar \cite[Appendix]{S-P} mathematically rigorously
justified the homogenization procedure. Using the same method of
two-scale expansions J. Keller and R. Burridge \cite{B-K} derived
formally the system of Biot's equations from the problem
(\ref{0.1})-- (\ref{0.7}) in the case when the   parameter $\alpha
_{\mu}$ was of order $\varepsilon^2$, and the rest of the
coefficients were fixed independent of $\varepsilon$.
 Under the same assumptions as in the article \cite{B-K}, the rigorous
justification of Biot's model was given by G. Nguetseng \cite{GNG}
and later by  Th. Clopeaut \emph{et al}. \cite{G-M2}.  The most
general case of the problem (\ref{0.1})-- (\ref{0.7}) when
$$\tau_{0} ,\, \mu_{0} , \, \lambda_{0}^{-1}, \, \nu_0, \,
p_{*}^{-1} , \, \eta_0^{-1}<\infty $$ has been studied in \cite{AM}.

All these authors have used Nguetseng's two-scale convergence method
\cite{NGU,LNW}.

In the present work by means of the same  method we investigate the
rest of all possible limiting regimes in the problem (\ref{0.1})--
(\ref{0.7}).  Namely, if $\tau_{0}=\infty $, which is a case of
short-time processes, then re-normalizing the displacement vector by
setting
$$ \w  \rightarrow  \alpha_\tau \w,$$
 we reduce the problem
to the case $\tau_{0}=1 $ and $\mu_{0}< \infty $. Here the only one
case $\lambda_{0}=0 $ needs an additional consideration.

 Therefore we restrict ourself by the case, when
 $$\nu_0, \, \mu_0< \infty ; \quad \lambda_0=0, \quad \tau_{0}=1, \quad 0< p_{*},\,\eta_0. $$

 We show that in the case $\mu_0>0$ the
homogenized equations are non-isotropic Stokes equations for fluid
velocity coupled with acoustic equations for the solid component, or
non-isotropic Stokes system for the one-velocity continuum (theorem
\ref{theorem2}). In the case $\mu_0=0$ the homogenized equations are
different types of acoustic equations for two-velocity or
one-velocity continuum (theorem \ref{theorem3}).

\section{Main results}

As usual, equation (\ref{0.1}) is understood in the sense of
distributions. It involves the proper equation (\ref{0.1}) in a
usual sense in the domains $\Omega_f^{\varepsilon}$ and
$\Omega_s^{\varepsilon}$ and the boundary conditions
\begin{eqnarray} \label{1.1}
& [\w]=0,\quad \x_0\in \Gamma ^{\varepsilon},\; t\geq 0,\\
\label{1.2} & [\textbf{P}]=0, \quad \x_0\in \Gamma ^{\varepsilon},\;
t\geq 0
\end{eqnarray}
on the boundary  $\Gamma^\varepsilon $, where
\begin{eqnarray}
\nonumber & [\varphi](\x_0)=\varphi_{(s)}(\x_0)
-\varphi_{(f)}(\x_0),\\
 \nonumber \displaystyle
& \varphi_{(s)}(\x_0) =\lim\limits_{\tiny \begin{array}{l}\x\to \x_0\\
\x\in \Omega_s^{\varepsilon}\end{array}} \varphi(\x),\quad
\varphi_{(f)}(\x_0) =\lim\limits_{\tiny \begin{array}{l}\x\to \x_0\\
\x\in \Omega_f^{\varepsilon}\end{array}} \varphi(\x).
\end{eqnarray}

 There are various equivalent in the sense of
distributions forms of representation of equation (\ref{0.1}) and
boundary conditions (\ref{1.1})--(\ref{1.2}). In what follows, it is
convenient to write them in the form of the integral equalities.

 We say that functions
$(\w^{\varepsilon},p^{\varepsilon},q^{\varepsilon},\pi^{\varepsilon})$
are called a generalized solution of the problem (\ref{0.1})--
(\ref{0.7}), if they satisfy the regularity conditions
\begin{equation} \label{1.3}
\w^{\varepsilon},\, D(x,\w^{\varepsilon}),\,
\div_x\w^{\varepsilon},\, q^{\varepsilon},\,p^{\varepsilon},\,
\frac{\partial p^{\varepsilon}}{\partial
 t},\, \pi ^{\varepsilon} \in L^2(\Omega_{T})
\end{equation}
in the domain $\Omega_{T}=\Omega\times (0,T)$, boundary conditions
(\ref{0.7}) in the trace sense, equations (\ref{0.3})-- (\ref{0.5})
a.e. in  $\Omega_{T}$  and  integral identity
\begin{eqnarray}\label{1.4}
\left. \begin{array}{lll} \displaystyle\int_{\Omega_{T}}
\Bigl(\alpha_\tau \rho ^{\varepsilon} \w^{\varepsilon}\cdot
\frac{\partial ^{2}{\mathbf \varphi}}{\partial t^{2}} - \chi
^{\varepsilon}\alpha_\mu D(\x, \w^{\varepsilon}):
D(x,\frac{\partial {\mathbf \varphi}}{\partial t})+\\[1ex]
\{(1-\chi^{\varepsilon})\alpha_\lambda D(x,\w^{\varepsilon})-
(q^{\varepsilon}+\pi^{\varepsilon})I\} : D(x,{\mathbf
\varphi})\Bigr) d\x dt=0
\end{array} \right\}
\end{eqnarray}
for all smooth  vector-functions  ${\mathbf \varphi}={\mathbf
\varphi}(\x,t)$ such that  ${\mathbf \varphi}|_{\partial \Omega}
={\mathbf \varphi}|_{t=T}=\partial {\mathbf \varphi} / \partial
t|_{t=T}=0$.

In (\ref{1.4})  by $A:B$ we denote the convolution (or,
equivalently, the inner tensor product) of two second-rank tensors
along the both indexes, i.e., $A:B=\mbox{tr\,} (B^*\circ
A)=\sum_{i,j=1}^3 A_{ij} B_{ji}$.

In what follows all parameters may take all permitted values.  If,
for example, $\eta_{0}^{-1}=0$ , then all terms in final equations
 containing this  parameter  disappear.

 The following theorems \ref{theorem1}--\ref{theorem3} are the main results of the paper.

\begin{theorem} \label{theorem1}
Let  $\F$, $\partial \F / \partial t $ and $\partial ^{2}\F /
\partial t^{2} $   are   bounded in $L^2(\Omega)$.  Then for all $\varepsilon >0$
  on the arbitrary time interval
$[0,T]$ there exists a unique generalized solution of the problem
(\ref{0.1})-- (\ref{0.7})  and
\begin{equation} \label{1.5}
 \displaystyle \max\limits_{0\leq t\leq
T}\| \frac{\partial ^{2}\w^{\varepsilon}}{\partial t^{2}}(t)
\|_{2,\Omega} \leq C_{0} ,
\end{equation}
\begin{equation} \label{1.6}
 \displaystyle \max\limits_{0\leq t\leq
T}\|\sqrt{\alpha_\mu} \chi^\varepsilon |\nabla_x \frac{\partial
\w^{\varepsilon}}{\partial t}(t)|+(1-\chi^\varepsilon)
 \sqrt{\alpha_\lambda}|\nabla_x\frac{\partial \w^{\varepsilon}}{\partial
t}(t)| \|_{2,\Omega} \leq C_{0} ,
\end{equation}
\begin{equation}\label{1.7}
 \|q^{\varepsilon}\|_{2,\Omega_{T}} +
\|p^{\varepsilon}\|_{2,\Omega_{T}} + \frac{\alpha _{\nu}}{\alpha
_{p}}\|\frac{\partial p^{\varepsilon}}{\partial
t}\|_{2,\Omega_{T}}\leq C_{0}
\end{equation}
where $C_{0}$ does not depend on the small parameter $\varepsilon $.
\end{theorem}

\begin{theorem} \label{theorem2}
Assume that the hypotheses in theorem \ref{theorem1} hold, and
$\mu_{0}>0$. Then functions $\partial \w^{\varepsilon} /
\partial t$  admit an extension $\vv^{\varepsilon}$ from
$\Omega_{f,T}^{\varepsilon}=\Omega_f^\varepsilon \times (0,T)$
 into $\Omega_{T}$ such that the sequence $\{\vv^{\varepsilon}\}$
 converges strongly in $L^{2}(\Omega_{T})$ and weakly in
 $L^{2}((0,T);W^1_2(\Omega))$ to the function $\vv$. At the same time,
 sequences $\{\w^\varepsilon\}$,  $\{(1-\chi)\w^\varepsilon\}$,
 $\{p^{\varepsilon}\}$, $\{q^{\varepsilon}\}$, and
 $\{\pi^{\varepsilon}\}$ converge weakly in $L^{2}(\Omega_{T})$
 to $\w$, $\w^{s}$,  $p$, $q$, and $\pi$, respectively.

\mathbf{I}) If $\lambda _{1}=\infty $, then $\partial\w^{s}/\partial
t=(1-m)\vv=(1-m)\partial\w/\partial t$ and weak and strong limits
$q$, $p$, $\pi$ and $\vv$ satisfy in $\Omega_{T}$ the
initial-boundary value problem
 \begin{equation}\label{1.8}
\left. \begin{array}{lll} \displaystyle \hat{\rho}\frac{\partial
\vv}{\partial t}=&& \div_x \{\mu
_{0}A^{f}_{0}:D(x,\vv) +  B^{f}_{0}\pi +B^{f}_{1}\div_x \vv+\\[1ex]
&&\int_{0}^{t}B^{f}_{2}(t-\tau)\div_x
\vv(\x,\tau)d\tau\}-\nabla(q+\pi )+\hat{\rho}\F,
\end{array} \right\}
\end{equation}
\begin{equation}\label{1.9}
\left. \begin{array}{lll} \displaystyle &&p_{*}^{-1}\partial p
/\partial t+C^{f}_{0}:D(x,\vv)+
a^{f}_{0}\pi +(a^{f}_{1}+m)\div_x\vv \\[1ex]
&&+\int_{0}^{t}a^{f}_{2}(t-\tau)\div_x \vv(\x,\tau)d\tau=0,
\end{array} \right\}
\end{equation}
\begin{equation}\label{1.10}
 q=p +\frac{\nu_0}{p_{*}}\frac{\partial p}{\partial t},\quad
\frac{1}{p_{*}}\frac{\partial p}{\partial t}
+\frac{1}{\eta_{0}}\frac{\partial \pi}{\partial t}+\div_x\vv=0,
\end{equation}
where $\hat{\rho}=m \rho_{f} + (1-m)\rho_{s}$, $m=\int_{Y}\chi dy$
and the symmetric strictly positively defined constant fourth-rank
tensor $A^{f}_{0}$, matrices  $C^{f}_{0}, B^{f}_{0}$, $B^{f}_{1}$
and $B^{f}_{2}(t)$ and scalars $a^{f}_{0}$, $a^{f}_{1}$ and
$a^{f}_{2}(t)$ are defined below by formulas (\ref{4.28}),
(\ref{4.30}) - (\ref{4.31}).

Differential equations (\ref{1.8}) are endowed with homogeneous
initial and boundary conditions
 \begin{equation}\label{1.11}
 \vv(\x,0)=0,\quad \x\in \Omega,
 \quad \vv(\x,t)=0, \quad \x\in S, \quad t>0.
\end{equation}

\mathbf{II}) If $\lambda _{1}<\infty $, then  weak and strong limits
$\w^{s}$, $q$, $p$, $\pi$ and $\vv$ satisfy in $\Omega_{T}$ the
initial-boundary value problem, which consists of Stokes like system
\begin{eqnarray}\label{1.12}
 \left. \begin{array}{lll} \displaystyle &&\rho_{f}m\partial
\vv / \partial t+\rho_{s}\partial ^2\w^{s} / \partial t^2 + \nabla
(q+\pi )-\hat{\rho}\F=\div_x \{B^{f}_{0}\pi+\\[1ex]
&&\mu_{0}A^{f}_{0}:D(x,\vv)+B^{f}_{1}\div_x \vv
+\int_{0}^{t}B^{f}_{2}(t-\tau)\div_x \vv(\x,\tau)d\tau\},
\end{array} \right\}
\end{eqnarray}
\begin{equation}\label{1.13}
\left. \begin{array}{lll} \displaystyle &&p_{*}^{-1}\partial p
/\partial t+C^{f}_{0}:D(x,\vv)+
a^{f}_{0}\pi +(a^{f}_{1}+m)\div_x\vv \\[1ex]
&&+\int_{0}^{t}a^{f}_{2}(t-\tau)\div_x \vv(\x,\tau)d\tau=0,
\end{array} \right\}
\end{equation}
\begin{equation}\label{1.14}
 q=p +\frac{\nu_0}{p_{*}}\frac{\partial p}{\partial t},
\end{equation}
for the liquid component coupled with a continuity equation
\begin{equation}\label{1.15}
\frac{1}{p_{*}}\frac{\partial p}{\partial t}
+\frac{1}{\eta_{0}}\frac{\partial \pi}{\partial t}+\div_x
\frac{\partial\w^{s}}{\partial t} +m\div_x \vv=0,
\end{equation}
 the relation
\begin{equation}\label{1.16}
\frac{\partial \w^{s}}{\partial t}=(1-m)\vv(\x,t)+\int_{0}^{t}
B^{s}_{1}(t-\tau)\cdot \z(\x,\tau )d\tau,
\end{equation}
$$\z(\x,t)=-\frac{1}{1-m}\nabla_x\pi(\x,t)+\rho_{s}\F(\x,t)-
\rho_{s}\frac{\partial \vv}{\partial t}(\x,t)$$ in the case of
$\lambda_{1}>0$, or the balance of momentum equation in the form
\begin{equation}\label{1.17}
\rho_{s}\frac{\partial^{2}\w^{s}}{\partial
t^{2}}=\rho_{s}B^{s}_{2}\cdot \frac{\partial \vv}{\partial
t}+((1-m)I-B^{s}_{2})\cdot(-\frac{1}{1-m}\nabla_x\pi+\rho_{s}\F)
\end{equation}
in the case of $\lambda_{1}=0$ for the solid component. The problem
is supplemented by boundary and initial conditions (\ref{1.11}) for
the velocity $\vv$ of the liquid component and by the homogeneous
initial conditions and the boundary condition
\begin{equation}\label{1.18}
\w^{s}(\x,t)\cdot \n(\x)=0,
     \quad (\x,t) \in S, \quad t>0,
\end{equation}
for the displacement $\w^{s}$ of the solid component. In Eqs.
(\ref{1.16})--(\ref{1.18}) $\n(\x)$ is the unit normal vector to $S$
at a point $\x \in S$, and matrices $B^{s}_{1}(t)$ and $B^{s}_{2}$
are given below by Eqs. (\ref{4.38}) and (\ref{4.40}), where the
matrix $((1-m)I - B^{s}_{2})$ is symmetric and positively definite.
\end{theorem}
\begin{theorem}\label{theorem3}
Assume that the hypotheses in Theorem \ref{theorem1} hold, and
$$\mu_{0}=0;\,\,\, \,  p_{*}, \, \eta_{0}<\infty.$$
Then  there exist  functions
$\w_{f}^\varepsilon,\,\w_{s}^\varepsilon \in L^\infty
(0,T;W^1_2(\Omega))$ such that
$$\w_{f}^\varepsilon = \w^\varepsilon \,\,
 \mbox{in}\,\,\Omega^{\varepsilon}_{f}\times (0,T), \quad \w_{s}^\varepsilon = \w^\varepsilon
  \,\, \mbox{in}\,\,\Omega^{\varepsilon}_{s}\times (0,T)$$
  and sequences $\{p^\varepsilon\}$, $\{q^\varepsilon\}$, $\{\pi^\varepsilon\}$, $\{\w^\varepsilon \}$,
 $\{\chi^{\varepsilon}\w^\varepsilon \}$, $\{\w_{f}^\varepsilon \}$ and $\{\w_{s}^\varepsilon \}$ converge
weakly in $L^2(\Omega_T)$  to functions  $p$, $q$, $\pi$, $\w$,
$\w^{f}$,  $\w_{f}$ and $\w_{s}$ respectively as
$\varepsilon\searrow 0$.

\textbf{I)}  If $\mu_{1}=\lambda_{1}=\infty$, then
$\w_{f}=\w_{s}=\w$ and functions $\w$, $p$, $q$   and $\pi$ satisfy
in $\Omega_{T}$ the system of acoustic equations
\begin{equation} \label{1.19}
\hat{\rho}\frac{\partial ^2\w}{\partial
t^2}=-\frac{1}{(1-m)}\nabla_x\pi+\hat{\rho}\F,
\end{equation}
\begin{equation} \label{1.20}
\frac{1}{p_{*}}p+\frac{1}{\eta_{0}}\pi +\div_x\w=0,
\end{equation}
\begin{equation} \label{1.21}
q=p+\frac{\nu_0}{p_{*}}\frac{\partial p}{\partial t}, \quad
\frac{1}{m}q=\frac{1}{1-m}\pi ,
\end{equation}
homogeneous initial conditions
\begin{equation} \label{1.22}
\w(\x,0)=\frac{\partial \w}{\partial t}(\x,0)=0, \quad \x \in \Omega
\end{equation}
and homogeneous boundary condition
\begin{equation} \label{1.23}
\w(\x,t)\cdot \n(\x)=0, \quad \x\in S,\, t>0.
\end{equation}
\textbf{II)}  If $\mu_{1}=\infty$ and $\lambda_{1}<\infty$, then
functions $ \w_{f}=\w$, $ \w^{s}$, $p$, $q$  and $\pi$  satisfy in
$\Omega_{T}$ the system of acoustic equations, which consist of the
state equations (\ref{1.21}) and  balance of momentum equation
\begin{equation} \label{1.24}
\rho_{f}m\frac{\partial ^2\w_{f}}{\partial
t^2}+\rho_{s}\frac{\partial ^2\w^{s}}{\partial
t^2}=-\frac{1}{(1-m)}\nabla_x\pi+\hat{\rho}\F,
\end{equation}
for the liquid component, continuity equation
\begin{equation} \label{1.25}
\frac{1}{p_{*}}p+\frac{1}{\eta_{0}}\pi
+m\div_x\w_{f}+\div_x\w^{s}=0,
\end{equation}
 and the relation
\begin{equation}\label{1.26}
\frac{\partial \w^{s}}{\partial t}=(1-m)\frac{\partial
\w_{f}}{\partial t}+\int_{0}^{t} B^{s}_{1}(t-\tau)\cdot
\z^{s}(\x,\tau )d\tau ,
\end{equation}
$$\z^{s}(\x,t)=-\frac{1}{1-m}\nabla_x\pi(\x,t)+\rho_{s}\F(\x,t)-
\rho_{s}\frac{\partial^{2}\w_{f}}{\partial t^{2}}(\x,t)$$ in the
case of $\lambda_{1}>0$, or the balance of momentum equation
 in the form
\begin{equation}\label{1.27}
\rho_{s}\frac{\partial^{2}\w^{s}}{\partial
t^{2}}=\rho_{s}B^{s}_{2}\cdot \frac{\partial^{2}\w_{f}}{\partial
t^{2}}+((1-m)I-B^{s}_{2})\cdot(-\frac{1}{1-m}\nabla_x\pi+\rho_{s}\F)
\end{equation}
in the case of $\lambda_{1}=0$ for the solid component. The problem
(\ref{1.21}), (\ref{1.24})--(\ref{1.27}) is supplemented by
homogeneous initial conditions (\ref{1.22}) for the displacements in
the liquid and the solid components and homogeneous boundary
condition (\ref{1.23}) for the displacements $\w=m\w_{f}+\w^{s}$.

In Eqs.(\ref{1.26})--(\ref{1.27}) matrices $B^{s}_{1}(t)$ and
$B^{s}_{2}$ are the same as in theorem \ref{theorem2}.

\textbf{III)}  If $\mu_{1}<\infty$ and $\lambda_{1}=\infty$, then
functions $ \w^{f}$, $ \w_{s}=\w$, $p$, $q$ and $\pi$  satisfy in
$\Omega_{T}$ the system of acoustic equations, which consist of the
state equations (\ref{1.21}) and the balance of  momentum equation
\begin{equation} \label{1.28}
\rho_{f}\frac{\partial ^2\w^{f}}{\partial
t^2}+\rho_{s}(1-m)\frac{\partial ^2\w_{s}}{\partial
t^2}=-\frac{1}{(1-m)}\nabla_x\pi+\hat{\rho}\F,
\end{equation}
for the solid component, the continuity equation
\begin{equation} \label{1.29}
\frac{1}{p_{*}}p+\frac{1}{\eta_{0}}\pi
+\div_x\w^{f}+(1-m)\div_x\w_{s}=0,
\end{equation}
 and the relation
\begin{equation}\label{1.30}
\frac{\partial \w^{f}}{\partial t}=m\frac{\partial \w_{s}}{\partial
t}+\int_{0}^{t} B^{f}_{1}(t-\tau)\cdot \z^{f}(\x,\tau )d\tau ,
\end{equation}
$$\z^{f}(\x,t)=-\frac{1}{m}\nabla_x q(\x,t)+\rho_{f}\F(\x,t)-
\rho_{f}\frac{\partial^{2}\w_{s}}{\partial t^{2}}(\x,t)$$ in the
case of $\lambda_{1}>0$, or the balance of momentum equation  in the
form
\begin{equation}\label{1.31}
\rho_{f}\frac{\partial^{2}\w^{f}}{\partial
t^{2}}=\rho_{f}B^{f}_{2}\cdot \frac{\partial^{2}\w_{s}}{\partial
t^{2}}+(mI-B^{f}_{2})\cdot(-\frac{1}{m}\nabla_x q+\rho_{f}\F)
\end{equation}
in the case of $\lambda_{1}=0$ for the liquid component. The problem
(\ref{1.21}), (\ref{1.28})--(\ref{1.31}) is supplemented by
homogeneous initial conditions (\ref{1.22}) for the displacements in
the liquid and the solid components and homogeneous boundary
condition (\ref{1.23}) for the displacements
$\w=\w^{f}+(1-m)\w_{s}$.

In Eqs.(\ref{1.30})--(\ref{1.31}) matrices $B^{f}_{1}(t)$ and
$B^{f}_{2}$ are given below by formulas (\ref{5.32})--(\ref{5.33}),
where the matrix $(mI - B^{f}_{2})$ is symmetric and positively
definite.

\textbf{IV)}  If $\mu_{1}<\infty$ and $\lambda_{1}<\infty$,  then
functions $ \w$, $p$, $q$  and $\pi$  satisfy in $\Omega_{T}$ the
system of acoustic equations, which consist of the continuity and
the state equations (\ref{1.20})  and (\ref{1.21}) and the relation
\begin{equation}\label{1.32}
\frac{\partial\w}{\partial
t}=\int_{0}^{t}B^{\pi}(t-\tau)\cdot\nabla\pi (\x,\tau )d\tau
+\textbf{f}(\x,t),
\end{equation}
where  $B^{\pi}(t)$ and $\textbf{f}(\x,t)$ are given below by Eqs.(
\ref{5.40}) and (\ref{5.41}).

The problem (\ref{1.20}), (\ref{1.21}), (\ref{1.32}) is supplemented
by homogeneous initial and boundary conditions (\ref{1.22}) and
(\ref{1.23}).
\end{theorem}

\section{Preliminaries}
\subsection{ Two-scale convergence}
Justification of theorems \ref{theorem1}--\ref{theorem3} relies on
systematic use of the method of two-scale convergence, which had
been proposed by G. Nguetseng \cite{NGU} and has been applied
recently to a wide range of homogenization problems (see, for
example, the survey \cite{LNW}).

\begin{definition} \label{TS}
A sequence $\{\varphi^\varepsilon\}\subset L^2(\Omega_{T})$ is said
to be \textit{two-scale convergent} to a limit $\varphi\in
L^2(\Omega_{T}\times Y)$ if and only if for any 1-periodic in $\y$
function $\sigma=\sigma(\x,t,\y)$ the limiting relation
\begin{equation}\label{2.1}
\lim_{\varepsilon\searrow 0} \int_{\Omega_{T}}
\varphi^\varepsilon(\x,t) \sigma\left(\x,t,\x /
\varepsilon\right)d\x dt = \int _{\Omega_{T}}\int_Y
\varphi(\x,t,\y)\sigma(\x,t,\y)d\y d\x dt
\end{equation}
holds.
\end{definition}

Existence and main properties of weakly convergent sequences are
established by the following fundamental theorem \cite{NGU,LNW}:
\begin{theorem} \label{theorem3.1}(\textbf{Nguetseng's theorem})

\textbf{1.} Any bounded in $L^2(Q)$ sequence contains a subsequence,
two-scale convergent to some limit
$\varphi\in L^2(\Omega_{T}\times Y)$.\\[1ex]
\textbf{2.} Let sequences $\{\varphi^\varepsilon\}$ and
$\{\varepsilon \nabla_x \varphi^\varepsilon\}$ be uniformly bounded
in $L^2(\Omega_{T})$. Then there exist a 1-periodic in $\y$ function
$\varphi=\varphi(\x,t,\y)$ and a subsequence
$\{\varphi^\varepsilon\}$ such that $\varphi,\nabla_y \varphi\in
L^2(\Omega_{T}\times Y)$, and $\varphi^\varepsilon$ and $\varepsilon
\nabla_x \varphi^\varepsilon$ two-scale converge to $\varphi$ and
$\nabla_y \varphi$,
respectively.\\[1ex]
\textbf{3.} Let sequences $\{\varphi^\varepsilon\}$ and $\{\nabla_x
\varphi^\varepsilon\}$ be bounded in $L^2(Q)$. Then there exist
functions $\varphi\in L^2(\Omega_{T})$ and $\psi \in
L^2(\Omega_{T}\times Y)$ and a subsequence from
$\{\varphi^\varepsilon\}$ such that $\psi$ is 1-periodic in $\y$,
$\nabla_y \psi\in L^2(\Omega_{T}\times Y)$, and
$\varphi^\varepsilon$ and $\nabla_x \varphi^\varepsilon$ two-scale
converge to $\varphi$ and $\nabla_x \varphi(\x,t)+\nabla_y
\psi(\x,t,\y)$, respectively.
\end{theorem}

\begin{corollary} \label{corollary3.1}
Let $\sigma\in L^2(Y)$ and
$\sigma^\varepsilon(\x):=\sigma(\x/\varepsilon)$. Assume that a
sequence $\{\varphi^\varepsilon\}\subset L^2(\Omega_{T})$ two-scale
converges to $\varphi \in L^2(\Omega_{T}\times Y)$. Then the
sequence $\sigma^\varepsilon \varphi^\varepsilon$ two-scale
converges to $\sigma \varphi$.
\end{corollary}
\subsection{An extension lemma}
The typical difficulty in homogenization problems while passing to a
limit in Model $B^\varepsilon$ as $\varepsilon \searrow 0$ arises
because of the fact that the bounds on the gradient of displacement
$\nabla_x \w^\varepsilon$ may be distinct in liquid and rigid
phases. The classical approach in overcoming this difficulty
consists of constructing of extension to the whole $\Omega$ of the
displacement field defined merely on $\Omega_s$. The following lemma
is valid due to the well-known results from \cite{ACE,JKO}. We
formulate it in appropriate for us form:

\begin{lemma} \label{Lemma1}
Suppose that assumption \ref{assumption1} on geometry of periodic
structure holds, $ \psi^\varepsilon\in W^1_2(\Omega^\varepsilon_s)$
and $\psi^\varepsilon =0$ on $S_{s}^{\varepsilon}=\partial\Omega
^\varepsilon_s \cap
\partial \Omega$ in the trace sense.  Then there exists a function
$ \sigma^\varepsilon \in
 W^1_2(\Omega)$ such that its restriction on the sub-domain
$\Omega^\varepsilon_s$ coincide with $\psi^\varepsilon$, i.e.,
\begin{equation} \label{2.2}
(1-\chi^\varepsilon(\x))( \sigma^\varepsilon(\x) - \psi^\varepsilon
(\x))=0,\quad \x\in\Omega,
\end{equation}
and, moreover, the estimate
\begin{equation} \label{2.3}
\|\sigma^\varepsilon\|_{2,\Omega}\leq C\|
\psi^\varepsilon\|_{2,\Omega ^{\varepsilon}_{s}}  , \quad \|\nabla_x
\sigma^\varepsilon\|_{2,\Omega} \leq  C \|\nabla_x
 \psi^\varepsilon\|_{2,\Omega ^{\varepsilon}_{s}}
\end{equation}
hold true, where the constant $C$ depends only on geometry $Y$ and
does not depend on $\varepsilon$.
\end{lemma}
\subsection{ Friedrichs--Poincar\'{e}'s inequality in periodic
structure} The following lemma was proved by L. Tartar in
\cite[Appendix]{S-P}. It specifies Friedrichs--Poincar\'{e}'s
inequality for $\varepsilon$-periodic structure.
\begin{lemma} \label{F-P}
Suppose that assumptions on the geometry of $\Omega^\varepsilon_f$
hold true. Then for any function $\varphi\in
\stackrel{\!\!\circ}{W^1_2}(\Omega^\varepsilon_f)$ the inequality
\begin{equation} \label{(F-P)}
\int_{\Omega^\varepsilon_f} |\varphi|^2 d\x \leq C \varepsilon^2
\int_{\Omega^\varepsilon_f} |\nabla_x \varphi|^2 d\x
\end{equation}
holds true with some constant $C$, independent of $\varepsilon$.
\end{lemma}
\subsection{ Some notation}
Further we denote

 1) $$ \langle\Phi \rangle_{Y} =\int_Y \Phi  dy, \quad
 \langle\Phi \rangle_{Y_{f}} =\int_Y \chi \Phi  dy,
 \quad
 \langle\Phi \rangle_{Y_{s}} =\int_Y (1-\chi )\Phi  dy,$$
$$\langle\varphi  \rangle_{\Omega } =\int_{\Omega } \varphi  dx,
\quad
  \langle\varphi  \rangle_{\Omega_{T}} =\int_{\Omega_{T}} \varphi
  dxdt.$$
2) If $\textbf{a}$ and $\textbf{b}$ are two vectors then the matrix
$\textbf{a}\otimes \textbf{b}$ is defined by the formula
$$(\textbf{a}\otimes \textbf{b})\cdot
\textbf{c}=\textbf{a}(\textbf{b}\cdot \textbf{c})$$ for any vector
$\textbf{c}$.

3) If $B$ and $C$ are two matrices, then $B\otimes C$ is a
forth-rank tensor such that its convolution with any matrix $A$ is
defined by the formula
$$(B\otimes C):A=B (C:A)$$.

4) By $I^{ij}$ we denote the $3\times 3$-matrix with just one
non-vanishing entry, which is equal to one and stands in the $i$-th
row and the $j$-th column.

5) We  also  introduce
$$J^{ij}=\frac{1}{2}(I^{ij}+I^{ji})=\frac{1}{2} ({\mathbf e}_i
\otimes {\mathbf e}_j + {\mathbf e}_j \otimes {\mathbf e}_i),$$
where $({\mathbf e}_1, {\mathbf e}_2, {\mathbf e}_3)$ are the
standard Cartesian basis vectors.
\section{Proof of theorem \ref{theorem1}}
Estimates
 (\ref{1.5})-(\ref{1.6})  follow from
\begin{equation}\label{3.1}
\left. \begin{array}{lll} \displaystyle
\max\limits_{0<t<T}(\sqrt{\alpha_\eta}\| \div_x
 \partial\w^{\varepsilon}/\partial t(t) \|_{2,\Omega _s^{\varepsilon}}+
 \sqrt{\alpha_\lambda}\|\nabla_x \partial\w^{\varepsilon}/\partial t(t) \|_{2,\Omega
_s^{\varepsilon}}\\[1ex]
+ \sqrt{\alpha_\tau}\| \partial ^{2}\w^\varepsilon/\partial
t^{2}(t)\|_{2,\Omega}+\sqrt{\alpha _{p}} \|  \div_x
\partial\w^{\varepsilon}/\partial t(t)\|_{2,\Omega
_f^{\varepsilon}})\\[1ex]
+\sqrt{\alpha_\mu}\|\chi ^{\varepsilon}
 \nabla_x \partial ^{2}\w^\varepsilon/\partial t^{2} \|_{2,\Omega_T}+
 \sqrt{\alpha _{\nu}}\| \chi ^{\varepsilon} \div_x
\partial ^{2}\w^\varepsilon/\partial t^{2}\|_{2,\Omega _{T}}
\leq C_{0}/\sqrt{\alpha_\tau},
\end{array} \right\}
\end{equation}
where  $C_{0}$ is independent of  $\varepsilon$. Last estimates we
obtain if we differentiate  equation  for $\w^{\varepsilon}$ with
respect to time, multiply  by $\partial ^{2} \w^{\varepsilon} /
\partial t^{2}$ and  integrate by parts using continuity and state
equations (\ref{0.3})-- (\ref{0.5}).  The same estimates guarantee
the existence and uniqueness of the generalized solution for the
problem (\ref{0.1})-- (\ref{0.7}).

If $p_{*}+\eta_0<\infty $, then estimate  (\ref{1.7})  for pressures
follows from estimate (\ref{3.1}) and continuity and state equations
(\ref{0.3})-- (\ref{0.5}).

For the case $p_{*}+\eta_0=\infty $ estimate  (\ref{1.7})  follows
from integral identity (\ref{1.4}) and estimates (\ref{3.1}) as an
estimate of the corresponding functional, if we re-normalized
pressures, such that
 $$\int _{\Omega} (q^\varepsilon(\x,t)+\pi^\varepsilon(\x,t)) d\x=0. $$.

Indeed, integral identity  (\ref{1.4}) and estimates  (\ref{3.1})
imply
$$|\int _{\Omega} (q^\varepsilon+\pi^\varepsilon )\div_x {\mathbf{\psi}} d\x
|\leq C \|\nabla {\mathbf{\psi}}\|_{2,\Omega}.$$

Choosing now ${\mathbf{\psi}}$ such that
$(q^\varepsilon+\pi^\varepsilon )= \div_x {\mathbf{\psi}}$ we get
the desired estimate for the sum of pressures
$(q^\varepsilon+\pi^\varepsilon )$.  Such a choice is always
possible (see \cite{LAD}), if we put
$${\mathbf{\psi}}=\nabla \varphi + {\mathbf{\psi_{0}}}, \quad
\div_x {\mathbf{\psi_{0}}}=0, \quad \triangle
\varphi=q^\varepsilon+\pi^\varepsilon ,\quad \varphi |
_{\partial\Omega}=0, \quad (\nabla \varphi + {\mathbf{\psi_{0}}})|
_{\partial\Omega}=0.$$
 Note that the re-normalization of the pressures  $(q^\varepsilon+\pi^\varepsilon )$
  transforms continuity  equations
  (\ref{0.4})-(\ref{0.5})  for  pressures into
\begin{eqnarray} \label{3.2}
&\displaystyle \frac{1}{\alpha_p}p^{\varepsilon}+
\chi^{\varepsilon}\div_x \w^{\varepsilon}=\frac{1}{m}\beta ^{\varepsilon}\chi^\varepsilon ,\\
\label{3.3}& \displaystyle \frac{1}{\alpha_\eta}\pi^{\varepsilon}
+(1-\chi^{\varepsilon})\div_x \w^{\varepsilon}=-\frac{1}{(1-m)}\beta
^{\varepsilon}(1-\chi^{\varepsilon}),
\end{eqnarray}
where
$$\beta ^{\varepsilon}=\langle \chi^\varepsilon \div_x
\w^\varepsilon \rangle _{\Omega}.$$

In what follows we will use equations (\ref{3.2}) and (\ref{3.3})
only if $p_{*}+\eta_0=\infty $.

Note that for the last case the basic integral identity (\ref{1.4})
permits to bound only the sum $(q^\varepsilon +\pi^{\varepsilon})$.
But thanks to the property that the product of these two functions
is equal to zero, it is enough to get bounds for each of these
functions. The pressure $p^{\varepsilon}$ is bounded from the  state
equation (\ref{0.3}), if we substitute the term  $(\alpha_{\nu} /
\alpha_p)\partial p^\varepsilon / \partial t$ from the continuity
equation \ref{3.2} and use estimate (\ref{3.1}).

\section{ Proof of theorem \ref{theorem2}}
\subsection{ Weak and two-scale limits of sequences of displacement
and pressures} On the strength of theorem \ref{theorem1}, the
sequences $\{p^\varepsilon\}$, $\{q^\varepsilon\}$,
$\{\pi^\varepsilon\}$  and  $\{\w^\varepsilon \}$   are uniformly in
$\varepsilon$ bounded in $L^2(\Omega_{T})$. Hence there exist a
subsequence of small parameters $\{\varepsilon>0\}$ and functions
$p$, $q$, $\pi$ and  $\w$  such that
\begin{equation}\label{4.1}
p^\varepsilon \rightarrow p,\quad q^\varepsilon \rightarrow q, \quad
\pi^\varepsilon \rightarrow \pi,  \quad  \w^\varepsilon \rightarrow
\w
\end{equation}
weakly in  $L^2(\Omega_T)$ as $\varepsilon\searrow 0$.

Moreover, due to lemma \ref{Lemma1} there is a function
$\vv^\varepsilon \in L^\infty (0,T;W^1_2(\Omega))$ such that
$\vv^\varepsilon =\partial \w^\varepsilon / \partial t$ in
$\Omega_{f}\times (0,T)$, and the family $\{\vv^\varepsilon \}$ is
uniformly in $\varepsilon$ bounded in $L^\infty
(0,T;W^1_2(\Omega))$. Therefore it is possible to extract a
subsequence of $\{\varepsilon>0\}$ such that
\begin{equation} \label{4.2}
\vv^\varepsilon \rightarrow \vv \mbox{ weakly in } L^2
(0,T;W^1_2(\Omega))
\end{equation}
as $\varepsilon \searrow 0$.

 Note also, that
\begin{equation} \label{4.3}
(1-\chi^\varepsilon )\alpha_\lambda D(x,\w^\varepsilon) \rightarrow
0.
\end{equation}
strongly in  $L^2(\Omega_T)$ as  $\varepsilon\searrow 0$.

 Relabeling if
necessary, we assume that the sequences converge themselves.

On the strength of Nguetseng's theorem, there exist 1-periodic in
$\y$ functions $P(\x,t,\y)$, $\Pi(\x,t,\y)$, $Q(\x,t,\y)$,
$\W(\x,t,\y)$  and $\V(\x,t,\y)$ such that the sequences
$\{p^\varepsilon\}$, $\{\pi^\varepsilon\}$, $\{q^\varepsilon\}$,
$\{\w^\varepsilon \}$  and $\{\nabla_x \vv^\varepsilon \}$ two-scale
converge to $P(\x,t,\y)$, $\Pi(\x,t,\y)$, $Q(\x,t,\y)$,
$\W(\x,t,\y)$  and $\nabla _{x}\vv +\nabla_{y}\V(\x,t,\y)$,
respectively.

Note that  the sequence  $\{\div_x \w^\varepsilon \}$ weakly
converges to $\div_x \w$ and $ \vv \in L^2
(0,T;\stackrel{\!\!\circ}{W^1_2}(\Omega)).$   Last assertion follows
from the Friedrichs--Poincar\'{e}'s inequality for $\vv^\varepsilon$
in the $\varepsilon$-layer  of the boundary $S$ and from convergence
of sequence $\{\vv^\varepsilon \}$  to $\vv$  strongly
 in $L^2(\Omega_{T})$ and weakly in $L^2 ((0,T);W^1_2(\Omega))$.\\
\subsection{Micro- and macroscopic equations I}

We start this section with macro- and microscopic equations,
connected with continuity equations.
\begin{lemma} \label{lemma4.1}
For all $ \x \in \Omega$ and $\y\in Y$ weak  and two-scale limits of
the sequences $\{p^\varepsilon\}$, $\{\pi^\varepsilon\}$,
$\{q^\varepsilon\}$, $\{\w^\varepsilon\}$, and $\{\vv^\varepsilon\}$
satisfy the relations
\begin{eqnarray} \label{4.4}
& \Pi=\pi (1-\chi) / (1-m);\\
\label{4.5} & q=p+\nu_0 p_{*}^{-1}\partial p /\partial t, \quad
Q=P+\nu_0 p_{*}^{-1}\partial P /\partial t;\\
 \label{4.6} & p_{*}^{-1}\partial p / \partial t+m\div_x\vv+
 \langle \div_y\V\rangle_{Y_{f}}=\partial \beta /\partial t;\\
\label{4.7} &p_{*}^{-1}\partial P/\partial t+\chi(\div_x\vv+
\div_y\V)=(\chi /m)\partial \beta /\partial t;\\
\label{4.8} &p/p_{*}+\pi /\eta_{0}+\div_x\w=0;\\
\label{4.9} & \w(\x,t)\cdot \n(\x)=0,
\quad \x\in S,\, t>0;\\
\label{4.10} & \div_y \W=0;\\
\label{4.11} & \partial \W/\partial t=\chi \vv+(1-\chi)\partial
\W/\partial t,
\end{eqnarray}
where $\beta =\langle\langle
\div_y\UU\rangle_{Y_{f}}\rangle_{\Omega}$, if
$p_{*}+\eta_{0}=\infty$  and $\beta=0$, if $p_{*}+\eta_{0}<\infty$
and  $\n(\x)$ is the unit normal vector to $S$ at a point $\x \in
S$.
\end{lemma}

\begin{proof}
In order to prove Eq.(\ref{4.4}), into Eq.(\ref{1.4}) insert a test
function ${\mathbf \psi}^\varepsilon =\varepsilon {\mathbf
\psi}\left(\x,t,\x / \varepsilon\right)$, where ${\mathbf
\psi}(\x,t,\y)$ is an arbitrary 1-periodic and finite on $Y_s$
function in $\y$. Passing to the limit as $\varepsilon \searrow 0$,
we get
\begin{equation} \label{4.12}
\nabla_y \Pi(\x,t,\y)=0, \quad \y\in Y_{s}.
\end{equation}

Next, fulfilling the two-scale limiting passage in equality
$$\chi^{\varepsilon}\pi^{\varepsilon}=0$$
we arrive at
$$\chi \Pi=0$$
which along with Eqs.(\ref{4.12}) justifies Eq. (\ref{4.4}).

Eqs. (\ref{4.5})--(\ref{4.9}) appear as the results of two-scale
limiting passages in Eqs. (\ref{0.3}), (\ref{3.2})-- (\ref{3.3})
with the proper test functions being involved. Thus, for example,
Eqs. (\ref{4.8}) and (\ref{4.9}) arise, if we consider the sum of
Eq.(\ref{3.2}) and Eq.(\ref{3.3})
\begin{equation}\label{4.13}
\frac{1}{\alpha_{p}}p^\varepsilon
+\frac{1}{\alpha_{\eta}}\pi^\varepsilon +\div_x \w^\varepsilon
=\frac{1}{m(1-m)}\beta ^{\varepsilon}(\chi^\varepsilon -m),
\end{equation}
multiply by an arbitrary function, independent of the ``fast''
variable $\x/\varepsilon$, and then pass to the limit as
$\varepsilon\searrow 0$. In order to prove Eq. (\ref{4.10}), it is
sufficient to consider the two-scale limiting relations in Eq.
(\ref{4.13}) as $\varepsilon \searrow 0$ with the test functions
$\varepsilon \psi \left(\x / \varepsilon\right) h(\x,t)$, where
$\psi$ and $h$ are arbitrary smooth test functions. In order to
prove Eq. (\ref{4.11}) it is sufficient to consider the two-scale
limiting relations in
$$\chi ^{\varepsilon}(\partial \w^{\varepsilon}/\partial
t-\vv^{\varepsilon})=0.$$
\end{proof}
\begin{corollary}\label{corollary4.1}
If $p_{*}+\eta_{0}=\infty$, then weak limits  $p$, $\pi$ and $q$
satisfy relations
\begin{equation}\label{4.100}
\langle p \rangle _{\Omega}=\langle \pi \rangle _{\Omega}=\langle q
\rangle _{\Omega}=0.
\end{equation}
\end{corollary}
\begin{lemma} \label{lemma4.2} For all $(\x,t) \in \Omega_{T}$
the relations
\begin{equation} \label{4.14}
\div_y \{\mu_0\chi (D(y,\V)+D(x,\vv))- (Q +\frac{(1-\chi)}{(1-m)}\pi
\cdot I\}=0,
\end{equation}
holds true.
\end{lemma}

\begin{proof}
Substituting a test function of the form ${\mathbf \psi}^\varepsilon
=\varepsilon {\mathbf \psi}\left(\x,t,\x / \varepsilon \right)$,
where ${\mathbf \psi}(\x,t,\y)$ is an arbitrary 1-periodic in $\y$
function vanishing on the boundary $S$, into integral identity
(\ref{1.4}), and passing to the limit as $\varepsilon \searrow 0$,
we arrive at  Eq. (\ref{4.14}).
\end{proof}

\begin{lemma} \label{lemma4.3}
Let $\hat{\rho}=m \rho_{f} + (1-m)\rho_{s}$. Then functions
$\w^{s}=\langle \W\rangle _{Y_{s}}$, $\vv$, $q$ and  $\pi$ satisfy
in $\Omega_{T}$ the system of macroscopic equations
\begin{eqnarray}\label{4.15}
    &&\rho_{f}m\frac{\partial \vv}{\partial t}+
    \rho_{s}\frac{\partial ^2\w^{s}}{\partial t^2}-\hat{\rho}\F=\\
    &&\div_x \{\mu _{0}(mD(x,\vv)+
    \langle D(y,\V)\rangle _{Y_{f}}
    )-(q+\pi )\cdot I \}.\nonumber
\end{eqnarray}
\end{lemma}
\begin{proof}
Eqs. (\ref{4.15}) arise as the limit of Eqs. (\ref{1.4}) with test
functions being finite in $\Omega_T$ and independent of
$\varepsilon$.
\end{proof}
\subsection{Micro- and macroscopic equations II}
In this section we derive macro- and microscopic equations for the
solid component.
\begin{lemma}\label{lemma4.14}
If $\lambda_{1}=\infty$, then the weak  limits of
$\{\vv^\varepsilon\}$ and $\{\partial \w^\varepsilon / \partial t\}$
coincide.
\end{lemma}
\begin{proof}
Let $\mathbf{\Psi}(\x,t,\y)$ be an arbitrary function periodic in
$\y$. The sequence $\{\beta^{\varepsilon}\}$, where
$$\beta^{\varepsilon}=\int_{\Omega}\sqrt{\alpha_{\lambda}}\nabla\w^\varepsilon
(\x,t)\mathbf{\Psi}(\x,t,\x /\varepsilon )dx,$$ is uniformly bounded
in $\varepsilon$. Therefore,
$$\int_{\Omega}\varepsilon \nabla\w^\varepsilon
\mathbf{\Psi}(\x,t,\x /\varepsilon
)dx=\frac{\varepsilon}{\sqrt{\alpha_{\lambda}}}\beta^{\varepsilon}\rightarrow
0$$
 as  $\varepsilon\searrow 0$, which is equivalent to
 $$\int_{\Omega}\int_{Y}\W(\x,t,\y)\nabla_{y}\mathbf{\Psi}(\x,t,\y)dxdy=0,$$
 or $\W(\x,t,\y)=\w(\x,t).$
\end{proof}

\begin{lemma} \label{lemma4.15}
Let $\lambda_1 <\infty$. Then the weak and two-scale limits $\pi$
and $\W$  satisfy the microscopic relations
\begin{equation}\label{4.16}
\rho_{s}\frac{\partial ^{2}\W}{\partial t^{2}}=
\lambda_{1}\triangle_y \W -\nabla_y R -\frac{1}{1-m}\nabla_x \pi
+\rho_{s}\F, \quad \y \in Y_{s},
\end{equation}
\begin{equation}\label{4.17}
    \frac{\partial \W}{\partial t}=\vv, \quad \y \in \gamma
\end{equation}
in the case $\lambda_{1}>0$, and relations
\begin{equation}\label{4.18}
    \rho_{s}\frac{\partial ^{2}\W}{\partial t^{2}}= -\nabla_y R
    -\frac{1}{1-m}\nabla _{x}\pi
    +\rho_{s}\F, \quad \y \in Y_{s},
\end{equation}
    \begin{equation}\label{4.19}
    ( \frac{\partial \W}{\partial t}-\vv)\cdot{\mathbf n}=0, \quad \y \in \gamma
\end{equation}
in the case $\lambda_{1}=0$.

In Eq. (\ref{4.19}) ${\mathbf n}$ is the unit normal to $\gamma$.
\end{lemma}

\begin{proof}
 Differential equations (\ref{4.16}) and (\ref{4.18}) follow
 as $\varepsilon\searrow 0$
 from integral equality (\ref{1.4}) with the test function ${\mathbf
\psi}={\mathbf \varphi}(x\varepsilon^{-1})\cdot h({\mathbf x},t)$,
where ${\mathbf \varphi}$ is solenoidal and finite in $Y_{s}$.

Boundary conditions (\ref{4.17}) are the consequences of the
two-scale convergence of $\{\alpha_{\lambda}^{\frac{1}{2}}\nabla_x
\w^{\varepsilon}\}$ to the function
$\lambda_{1}^{\frac{1}{2}}\nabla_y\W(\x,t,\y)$. On the strength of
this convergence, the function $\nabla_y \W (\x,t,\y)$ is
$L^2$-integrable in $Y$. The boundary conditions (\ref{4.19}) follow
from Eqs. (\ref{4.10})-(\ref{4.11}).
\end{proof}
\subsection{Homogenized equations I}

Here we derive homogenized equations for the liquid component.
 \begin{lemma} \label{lemma4.6} If $\lambda_1 =\infty$  then $\partial
\w / \partial t=\vv$ and  the weak limits $\vv$, $p$, $q$, and $\pi$
satisfy in $\Omega_{T}$ the initial-boundary value problem
 \begin{equation}\label{4.20}
\left. \begin{array}{lll} \displaystyle \hat{\rho}\frac{\partial
\vv}{\partial t}=&& \div_x \{\mu
_{0}A^{f}_{0}:D(x,\vv) +  B^{f}_{0}\pi +B^{f}_{1}\div_x \vv+\\[1ex]
&&\int_{0}^{t}B^{f}_{2}(t-\tau)\div_x
\vv(\x,\tau)d\tau\}-\nabla(q+\pi )+\hat{\rho}\F,
\end{array} \right\}
\end{equation}
\begin{equation}\label{4.21}
\left. \begin{array}{lll} \displaystyle &&p_{*}^{-1}\partial p
/\partial t+C^{f}_{0}:D(x,\vv)+
a^{f}_{0}\pi +(a^{f}_{1}+m)\div_x\vv \\[1ex]
&&+\int_{0}^{t}a^{f}_{2}(t-\tau)\div_x \vv(\x,\tau)d\tau=0,
\end{array} \right\}
\end{equation}
\begin{equation}\label{4.22}
 q=p +\frac{\nu_0}{p_{*}}\frac{\partial p}{\partial t},\quad
\frac{1}{p_{*}}\frac{\partial p}{\partial t}
+\frac{1}{\eta_{0}}\frac{\partial \pi}{\partial t}+\div_x\vv=0,
\end{equation}
where the symmetric strictly  positively defined constant
fourth-rank tensor $A^{f}_{0}$, matrices  $C^{f}_{0}, B^{f}_{0}$,
$B^{f}_{1}$ and $B^{f}_{2}(t)$ and scalars $a^{f}_{0}$, $a^{f}_{1}$
and $a^{f}_{2}(t)$ are defined below by formulas (\ref{4.28}),
(\ref{4.30}) - (\ref{4.31}).

Differential equations (\ref{4.20}) are endowed with homogeneous
initial and boundary conditions
 \begin{equation}\label{4.23}
 \vv(\x,0)=0,\quad \x\in \Omega,
 \quad \vv(\x,t)=0, \quad \x\in S, \quad t>0.
\end{equation}

\end{lemma}

\begin{proof}
In the first place let us notice that $\vv =\partial \w / \partial
t$ due to lemma \ref{lemma4.14}.

The homogenized equations (\ref{4.20}) follow from the macroscopic
equations (\ref{4.15}), after we insert in them the expression
$$\mu_{0}\langle D(y,\V)\rangle _{Y_{f}}=
\mu_{0}A^{f}_{1}:D(x,\vv) + B^{f}_{0}\pi +B^{f}_{1}\div_x \vv+
 \int_{0}^{t}B^{f}_{2}(t-\tau)\div_x \vv(\x,\tau)d\tau.$$
In turn, this expression follows by virtue of solutions of Eq.
 (\ref{4.5}) in the form
$$Q=P-\nu_0\chi(\div_x\vv+
\div_y\V)+\nu_0(\chi /m)\partial \beta /\partial t$$
  and Eqs.(\ref{4.6}) and (\ref{4.14}) on the pattern cell $Y_{f}$.
 Indeed, setting
$$  \V=\sum_{i,j=1}^{3}\V^{(ij)}(\y)D_{ij}
 +\V^{(0)}(\y)\pi+\V^{(1)}(\y)\div_x \vv$$
 $$+\int_{0}^{t}\V^{(2)}(\y,t-\tau)\div_x \vv(\x,\tau)d\tau,$$

$$ Q =\mu_{0}\sum_{i,j=1}^{3}Q^{(ij)}(\y)D_{ij}
 +Q^{(0)}(\y)\pi+Q^{(1)}(\y)\div_x \vv $$
$$+\int_{0}^{t}Q^{(2)}(\y,t-\tau)\div_x \vv(\x,\tau)d\tau,$$

 $$P =\int_{0}^{t}P^{(2)}(\y,t-\tau)\div_x \vv(\x,\tau)d\tau,$$
  where
 $$D_{ij}(\x,t)=\frac{1}{2}(\frac{\partial v_{i}}{\partial x_{j}}(\x,t)+
 \frac{\partial v_{j}}{\partial x_{i}}(\x,t)),$$
we arrive at the following periodic-boundary value problems in $Y$:
\begin{equation}\label{4.24}
\left. \begin{array}{lll}  \displaystyle \div_y\{\chi (
D(y,\V^{(ij)})+J^{ij} - Q^{(ij)}\cdot I)\}=0,\\[1ex]
\mu_{0}Q^{(ij)}+\nu_0\chi\div_y \V^{(ij)}=0 \quad \mbox{for}\quad p_{*}<\infty,\\[1ex]
\chi\div_y \V^{(ij)} =0 \quad \mbox{for} \quad p_{*}=\infty;
\end{array} \right\}
\end{equation}
\begin{equation}\label{4.25}
\left. \begin{array}{lll}  \displaystyle \div_y
\{\mu_{0}\chi D(y,\V^{(0)}) -( Q^{(0)}+\frac{1-\chi}{1-m})\cdot I \}=0,\\[1ex]
Q^{(0)}+\nu_0\chi\div_y \V^{(0)}=0 \quad \mbox{for}\quad p_{*}<\infty,\\[1ex]
\chi\div_y \V^{(0)}=0 \quad \mbox{for} \quad p_{*}=\infty;
\end{array} \right\}
\end{equation}
\begin{equation}\label{4.26}
\left. \begin{array}{lll}  \displaystyle \div_y
\{\mu_{0}\chi D(y,\V^{(1)}) - Q^{(1)}\cdot I \}=0,\\[1ex]
Q^{(1)}+\nu_0\chi(\div_y\V^{(1)}+1)=0 \quad \mbox{for}\quad p_{*}<\infty,\\[1ex]
\chi(\div_y\V^{(1)}+1)=0 \quad \mbox{for} \quad p_{*}=\infty;
\end{array} \right\}
\end{equation}
$\V^{(2)}=0, \,Q^{(2)}=P^{(2)}=0$ for $p_{*}=\infty$ and
\begin{equation}\label{4.27}
\left. \begin{array}{lll}  \displaystyle \div_y \{\mu_{0}\chi
D(y,\V^{(2)}) -Q^{(2)}\cdot I
\}=0,\\[1ex]
Q^{(2)}=P^{(2)}-\nu_0 \chi\div_y\V^{(2)}, \\[1ex]
1/p_{*}\partial P^{(2)} /\partial t +\chi\div_y \V^{(2)} =0,
\,P^{(2)}(\y,0)=p_{*}
\end{array} \right\}
\end{equation}
for $p_{*}<\infty$.

 Note that for  $p_{*}+\eta_{0}=\infty$
  $$\beta=\langle\langle \div_y \V\rangle_{\Omega}\rangle_{Y_{f}}=
  \sum_{i,j=1}^{3}\langle \div_y\V^{(ij)}\rangle_{Y_{f}}
  \langle D_{ij}\rangle_{\Omega} +$$
  $$\langle \div_y\V^{(0)}\rangle_{Y_{f}}
   \langle \pi\rangle_{\Omega}+ \langle \div_y\V^{(1)}\rangle_{Y_{f}}\langle \div_x \vv\rangle_{\Omega}
  =0$$
 due to homogeneous boundary conditions for  $\vv(\x,t)$ and Eq. (\ref{4.100}).

On the strength of the assumptions on the geometry of the pattern
``liquid'' cell $Y_{f}$, problems  (\ref{4.24})-- (\ref{4.26}) have
unique solution, up to an arbitrary constant vector. In order to
discard the arbitrary constant vectors we demand
$$
\langle\V^{(ij)}\rangle _{Y_{f}}=\langle\V^{(0)}\rangle_{Y_{f}}
=\langle\V^{(1)}\rangle_{Y_{f}}=\langle\V^{(2)}\rangle_{Y_{f}}=0.
$$
Thus
 \begin{equation}\label{4.28}
 A^{f}_{0}=\sum_{i,j=1}^{3}J^{ij}\otimes J^{ij} + A^{f}_{1}, \quad
 A^{f}_{1}=\sum_{i,j=1}^{3}\langle D(y,\V^{(ij)})\rangle _{Y_{f}}\otimes
    J^{ij}.
\end{equation}
 Symmetry of the tensor $A^{f}_{0}$ follows from symmetry of
 the tensor $A^{f}_{1}$. And symmetry of the latter one follows
 from the equality
 \begin{eqnarray}\label{4.29}
     &&\langle D(y,\V^{(ij)})\rangle _{Y_{f}} : J^{kl}
 =-\langle D(y,\V^{(ij)}) : D(y,\V^{(kl)})\rangle_{Y_{f}}\\
 && -\mbox{sgn}(\frac{1}{p_{*}})\frac{\nu_{0}}{\mu_{0}}\div_y
\V^{(ij)}\div_y \V^{(kl)},\nonumber
\end{eqnarray}
which appears by means of multiplication of Eq. (\ref{4.24}) for
$\V^{(ij)}$ by $\V^{(kl)}$ and by integration by parts using the
continuity equation.

This equality also implies positive definiteness of the tensor
$A^{f}_{0}$. Indeed, let $\zeta$ be an arbitrary symmetric matrix.
Setting
 $$\textbf{Z}=\sum_{i,j=1}^{3}\V^{(ij)}\zeta_{ij}$$
and taking into account Eq.(\ref{4.29}) we get

$$\langle D(y,\textbf{Z})\rangle _{Y_{f}}:\zeta
 =-\langle D(y,\textbf{Z}): D(y, \textbf{Z})\rangle_{Y_{f}} -
 \mbox{sgn}(\frac{1}{p_{*}})\frac{\nu_{0}}{\mu_{0}}(\div_y \textbf{Z})^{2},$$
\
 This equality and the definition of the tensor $A_{0}^f$ give

 $$(A_{0}^f:\zeta ):\zeta =
 \langle(D(y,\textbf{Z})+\zeta ): (D(y,\textbf{Z})+\zeta )\rangle_{Y_{s}}
 +\mbox{sgn}(\frac{1}{p_{*}})\frac{\nu_{0}}{\mu_{0}}(\div_y \textbf{Z})^{2}.$$

Now the strict positive definiteness of the tensor $A_{0}^{f}$
follows from the equality immediately above and the geometry of the
elementary cell $Y_{f}$. Namely, suppose that $(A_{0}^{s}:\zeta
):\zeta =0$ for some function $\zeta$, such that $\zeta :\zeta =1$.
Then  $(D(y,\textbf{Z})+\zeta )=0$, which is possible iff
$\textbf{Z}$ is a linear function in $\y$. On the other hand, all
linear periodic functions on $Y_{f}$ are constant. Finally, the
normalization condition $\langle\V^{(ij)}\rangle_{Y_{f}} =0$ yields
that $\textbf{Z}=0$. However, this is impossible because the
functions $\V^{(ij)}$ are linearly independent.

 Finally, Eqs. (\ref{4.21}) and (\ref{4.22}) for the pressures follow from Eqs.
 (\ref{4.5}), (\ref{4.6}), (\ref{4.8}) and
 $$\langle \div_y\V\rangle_{Y_{f}}=C^{f}_{0}:D(x,\vv)+
a^{f}_{0}\pi +a^{f}_{1}\div_x \vv+
\int_{0}^{t}a^{f}_{2}(t-\tau)\div_x \vv(\x,\tau)d\tau$$
 with
\begin{equation}\label{4.30}
  B^{f}_{i}=\mu_{0}\langle D(y,\V^{(i)})\rangle _{Y_{f}}, \quad
  i=0,1,2,
\end{equation}
\begin{equation}\label{4.31}
 C^{f}_{0}=\sum_{i,j=1}^{3}\langle \div_y\V^{(ij)}\rangle _{Y_{f}}J^{ij}, \quad
  a^{f}_{i}=\langle\div_y\V^{(i)}\rangle _{Y_{f}},\quad i=0,1,2.
\end{equation}
\end{proof}

\subsection{ Homogenized equations II}
We complete  the proof of theorem \ref{theorem2} with homogenized
equations for the solid component.

Let $\lambda_{1}<\infty$. In the same manner as above, we verify
that the limit $\vv$ of the sequence $\{\vv^\varepsilon\}$ satisfies
the initial-boundary value problem likes (\ref{4.20})--
(\ref{4.23}). The main difference here that, in general, the weak
limit $\partial\w / \partial t$ of the sequence
$\{\partial\w^\varepsilon /\partial t\}$ differs from $\vv$. More
precisely, the following statement is true.
\begin{lemma} \label{lemma4.7}
If $\lambda_{1}<\infty$ then the weak limits $\vv$, $\w^{s}$, $p$,
$q$, and $\pi$ of the sequences
 $\{\vv^\varepsilon\}$, $\{(1-\chi^{\varepsilon})\w^\varepsilon\}$,
  $\{p^\varepsilon\}$,  $\{q^\varepsilon\}$, and $\{\pi^\varepsilon\}$
satisfy the initial-boundary value problem in $\Omega_T$, consisting
of the balance of momentum equation
\begin{eqnarray}\label{4.32}
 &&\rho_{f}m\frac{\partial
\vv}{\partial t}+\rho_{s}
\frac{\partial ^2\w^{s}}{\partial t^2} + \nabla (q+\pi )-\hat{\rho}\F= \\
&&\div_x \{\mu_{0}A^{f}_{0}:D(x,\vv) + B^{f}_{0}\pi
 +B^{f}_{1}\div_x \vv \}+\int_{0}^{t}B^{f}_{2}(t-\tau)\div_x
\vv(\x,\tau)d\tau\},\nonumber
\end{eqnarray}
 the continuity equation (\ref{4.21}) and  first state equation
 in (\ref{4.22}) for the liquid component,
where $A^{f}_{0}$, $B^{f}_{0}$-- $B^{f}_{2}$ are the same as in
(\ref{4.20}),  the continuity equation
\begin{equation} \label{4.33}
\frac{1}{p_{*}}\frac{\partial p}{\partial t}+
\frac{1}{\eta_{0}}\frac{\partial\pi}{\partial t}+\div_x
\frac{\partial\w^{s}}{\partial t} +m\div_x \vv=0,
   \end{equation}
 the relation
\begin{equation}\label{4.34}
\frac{\partial \w^{s}}{\partial t}=(1-m)\vv(\x,t)+\int_{0}^{t}
B^{s}_{1}(t-\tau)\cdot \z(\x,\tau )d\tau ,
\end{equation}
$$\z(\x,t)=-\frac{1}{1-m}\nabla_x\pi(\x,t)+\rho_{s}\F(\x,t)-
\rho_{s}\frac{\partial \vv}{\partial t}(\x,t)$$ in the case of
$\lambda_{1}>0$, or the balance of momentum equation  in the form
\begin{equation}\label{4.35}
\rho_{s}\frac{\partial^{2}\w^{s}}{\partial
t^{2}}=\rho_{s}B^{s}_{2}\cdot \frac{\partial \vv}{\partial
t}+((1-m)I-B^{s}_{2})\cdot(-\frac{1}{1-m}\nabla_x\pi+\rho_{s}\F)
\end{equation}
in the case of $\lambda_{1}=0$ for the solid component. The problem
is supplemented by boundary and initial conditions (\ref{4.23}) for
the velocity $\vv$ of the liquid component and by homogeneous
initial conditions and the boundary condition
\begin{equation}\label{4.36}
\w^{s}(\x,t)\cdot \n(\x)=0, \quad (\x,t) \in S, \quad t>0,
\end{equation}
for the displacement $\w^{s}$ of the solid component. In Eqs.
(\ref{4.34})--(\ref{4.36}) $\n(\x)$ is the unit normal vector to $S$
at a point $\x \in S$, and matrices $B^{s}_{1}(t)$ and $B^{s}_{2}$
are given below by Eqs. (\ref{4.38}) and (\ref{4.40}).
\end{lemma}
\begin{proof}
The boundary condition (\ref{4.36}) follows from Eq.(\ref{4.9}),
equality
$$\frac{\partial\w}{\partial t}=\frac{\partial\w^{s}}{\partial t}+m\vv,$$
and homogeneous boundary condition for $\vv$.

 The same equality and Eq.(\ref{4.8}) imply (\ref{4.33}).
 The homogenized equations of balance of momentum (\ref{4.32})
  derives exactly as before.  Therefore we
omit the relevant proofs now and focus ourself only on derivation of
homogenized equation of the balance of momentum for the solid
displacements $\w^{s}$.

a) If  $\lambda_{1}>0$, then the solution of the system of
microscopic equations (\ref{4.10}), (\ref{4.16}), and (\ref{4.17}),
provided with the homogeneous initial data, is given by formula

$$\W=\int_{0}^{t}( \vv(\x,\tau )+
\textbf{B}^{s}_{1}(\y,t-\tau)\cdot\z)(\x,\tau )d\tau ,\quad
R=\int_{0}^{t}{\mathbf R}_{f}(\y,t-\tau)\cdot\z(\x,\tau )d\tau ,$$

 in which

$$\textbf{B}^{s}_{1}(\y,t)= \sum_{i=1}^{3}\W^{i}(\y,t)\otimes {\mathbf e}_{i},
\quad {\mathbf R}_{f}(\y,t)=\sum_{i=1}^{3}R^{i}(\y,t){\mathbf
e}_{i},$$

and the functions $\W^{i}(\y,t)$ and $R^{i}(\y,t)$  are defined by
virtue of the periodic initial-boundary value problem
\begin{equation}\label{4.37}
\left. \begin{array}{lll}  \displaystyle\rho_{s}\frac{\partial ^{2}
\W^{i}}{\partial t^{2}}-\lambda_{1}\triangle \W^{i} +\nabla R^{i}
=0,
\quad \mbox{div}_y \W^{i} =0, \quad \y \in Y_{s},\,  t>0,\\[1ex]
\W^{i}=0, \quad \y \in \gamma ,\,\,  t>0;\\[1ex]
 \W^{i}(y,0)=0, \quad
 \rho_{s}\partial \W^{i}/\partial t(y,0)={\mathbf e}_{i},\quad \y \in Y_{s}.
 \end{array} \right\}
\end{equation}
 In Eq. (\ref{4.37}) ${\mathbf e}_{i}$ is the standard Cartesian
 basis vector.

Therefore
\begin{equation}\label{4.38}
B^{s}_{1}(t)= \langle \frac{\partial \textbf{B}^{s}_{1}}{\partial
t}\rangle _{Y_{s}}(t).
\end{equation}

Note, that equations (\ref{4.37}) are understood in the sense of
distributions and the function $B^{s}_{1}(t)$ has no time derivative
at $t=0$.

 b) If  $\lambda_{1}=0$ then in the
process of solving the system (\ref{4.10}), (\ref{4.18}), and
(\ref{4.19}) we firstly find the pressure $R(\x,t,\y)$ by virtue of
solving the Neumann problem for Laplace's equation in
 $Y_{s}$ in the form

 $$R(\x,t,\y)=\sum_{i=1}^{3} R_{i}(\y) {\mathbf e}_{i}\cdot \z(\x,t),$$
 where $R^{i}(\y)$ is the solution of the problem
 \begin{equation}\label{4.39}
\triangle_y R_{i}=0,\quad \y \in Y_{s}; \quad \nabla_y R_{i}\cdot \n
=\n\cdot{\mathbf e}_{i}, \quad \y \in \gamma .
\end{equation}
Formula (\ref{4.35}) appears as the result of  homogenization of
Eq.(\ref{4.18}) and
 \begin{equation}\label{4.40}
B^{s}_{2}=\sum_{i=1}^{3}\langle \nabla R_{i}(\y)\rangle
_{Y_{s}}\otimes {\mathbf e}_{i},
\end{equation}
where the matrix $((1-m)I - B^{s}_{2})$ is symmetric and positively
definite.  In fact, let $\tilde{R}=\sum_{i=1}^{3}R_{i}\xi_{i}$ for
any unit vector $\xi$. Then
 $$(B\cdot\xi)\cdot\xi=\langle(\xi-\nabla\tilde{R})^{2}\rangle_{Y_{f}}>0$$
 due to the same reasons as in lemma \ref{lemma4.6}.

Note, that on the strength of the assumptions on the geometry of the
pattern ``solid'' cell $Y_{s}$, problem  (\ref{4.37}) has unique
solution and problem (\ref{4.39}) has unique solution up to an
arbitrary constant.
\end{proof}
\section{ Proof of theorem \ref{theorem3}}
\subsection{Weak and two-scale limits of sequences of displacement
and pressures} Let $\mu_{0}=0$. We again use Lemma \ref{Lemma1} and
conclude that there are  functions
$\w_{f}^\varepsilon,\,\w_{s}^\varepsilon \in L^\infty
(0,T;W^1_2(\Omega))$ such that
$$\w_{f}^\varepsilon = \w^\varepsilon \,\,
 \mbox{in}\,\,\Omega_{f}\times (0,T), \quad \w_{s}^\varepsilon = \w^\varepsilon
  \,\, \mbox{in}\,\,\Omega_{s}\times (0,T).$$

 On the strength of theorem \ref{theorem1}, the sequences $\{p^\varepsilon\}$,
$\{q^\varepsilon\}$, $\{\pi^\varepsilon\}$, $\{\w^\varepsilon \}$,
$\{\w_{f}^\varepsilon \}$,
$\{\sqrt{\alpha_{\mu}}\nabla\w_{f}^\varepsilon \}$,
$\{\w_{s}^\varepsilon \}$ and
 $\{\sqrt{\alpha_{\lambda}}\nabla\w_{s}^\varepsilon \}$  are
uniformly in $\varepsilon$ bounded in $L^2(\Omega_{T})$. Hence there
exist a subsequence of small parameters $\{\varepsilon>0\}$ and
functions $p$, $q$, $\pi$, $\w$, $\w_{f}$ and $\w_{s}$  such that
\begin{equation}\label{5.1}
p^\varepsilon \rightarrow p,\quad q^\varepsilon \rightarrow q, \quad
\pi^\varepsilon \rightarrow \pi,\quad  \w^\varepsilon \rightarrow
\w, \quad  \w_{f}^\varepsilon \rightarrow \w_{f},\quad
\w_{s}^\varepsilon \rightarrow \w_{s}
\end{equation}
weakly in  $L^2(\Omega_T)$ as $\varepsilon\searrow 0$.

 Note also, that
\begin{equation} \label{5.2}
(1-\chi^\varepsilon )\alpha_\lambda D(x,\w_{s}^\varepsilon)
\rightarrow 0, \quad \chi^\varepsilon \alpha_\mu
 D(x,\w_{f}^\varepsilon) \rightarrow 0
\end{equation}
strongly in  $L^2(\Omega_T)$ as  $\varepsilon\searrow 0$.

 Relabeling if necessary, we assume that the sequences converge themselves.

On the strength of Nguetseng's theorem, there exist 1-periodic in
$\y$ functions $P(\x,t,\y)$, $\Pi(\x,t,\y)$, $Q(\x,t,\y)$,
$\W(\x,t,\y)$,  $\W_{f}(\x,t,\y)$  and $\W_{s}(\x,t,\y)$ such that
the sequences $\{p^\varepsilon\}$, $\{\pi^\varepsilon\}$,
$\{q^\varepsilon\}$, $\{\w^\varepsilon \}$, $\{\w_{f}^\varepsilon
\}$  and $\{\w_{s}^\varepsilon \}$ two-scale converge to
$P(\x,t,\y)$, $\Pi(\x,t,\y)$, $Q(\x,t,\y)$, $\W(\x,t,\y)$,
$\W_{f}(\x,t,\y)$  and $\W_{s}(\x,t,\y)$, respectively.

 Finally note, that if $\mu_{1}=\infty$ ($\lambda_{1}=\infty$), then due to Lemma
\ref{lemma4.14} the sequence $\{\w_{f}^\varepsilon \}$
($\{\w_{s}^\varepsilon \}$) converges strongly to $\w_{f}$
($\w_{s}$) and $\w^{f}=\langle\W\rangle_{Y_{f}}=m\w_{f}$ (
$\w^{s}=\langle\W\rangle_{Y_{s}}=(1-m)\w_{s}$).\\

\subsection{Micro- and macroscopic equations}
As before we start the proof of theorem with  macro- and microscopic
equations, connected with continuity equations.
\begin{lemma} \label{lemma5.1}
For all $ \x \in \Omega$ and $\y\in Y$ weak  and two-scale limits of
the sequences $\{p^\varepsilon\}$, $\{\pi^\varepsilon\}$,
$\{q^\varepsilon\}$, $\{\w^\varepsilon\}$,
 $\{\w_{f}^\varepsilon\}$, and $\{\w_{s}^\varepsilon\}$ satisfy the
relations
\begin{eqnarray} \label{5.3}
& Q=q\chi /m, \quad P=p\chi /m, \quad \Pi=\pi (1-\chi) / (1-m);\\
\label{5.4} & q/m=\pi /(1-m), \quad
q=p+\nu_0 p_{*}^{-1}\partial p /\partial t;\\
\label{5.5} & p/p_{*}+\pi /\eta_{0}+\div_x\w=0;\\
\label{5.6} & \w(\x,t)\cdot \n(\x)=0,
\quad \x\in S,\, t>0;\\
\label{5.7} & \div_y \W=0;\\
\label{5.8} & \W=\chi\W_{f}+(1-\chi)\W_{s}.
\end{eqnarray}
\end{lemma}

\begin{proof}
The derivation of  Eqs.(\ref{5.3})--(\ref{5.8}) is the same as the
derivation of Eqs.(\ref{4.4}), (\ref{4.5}), (\ref{4.8}) and
(\ref{4.10})--(\ref{4.11}) in lemma \ref{lemma4.1}. Thus, for
example, the first relation in Eq.(\ref{5.4}) follow from
Eq.(\ref{5.3}) and from the strong convergence of the sequence
$\{(q^{\varepsilon}+\pi^{\varepsilon})$ to $(q+\pi)$, which implies
the equality $Q+\Pi=q+\pi$.
\end{proof}

\begin{lemma} \label{lemma5.2} For all $(\x,t) \in \Omega_{T}$ the relation
\begin{equation} \label{5.9}
\rho_{f}\frac{\partial ^2\w^{f}}{\partial
t^2}+\rho_{s}\frac{\partial ^2\w^{s}}{\partial
t^2}=-\frac{1}{(1-m)}\nabla_x\pi+\hat{\rho}\F,
\end{equation}
 holds true.
\end{lemma}

\begin{proof}
Substituting a test function of the form ${\mathbf \psi} ={\mathbf
\psi}(\x,t)$ into integral identity (\ref{1.4}), and passing to the
limit as $\varepsilon \searrow 0$, we arrive at Eq.(\ref{5.9}).
\end{proof}

\begin{lemma} \label{lemma5.3}
Let $\mu_{1}=\infty$ and $\lambda_{1}<\infty$. Then functions $\{\W,
\w_{f}, \,  \pi \}$ satisfy in $Y_{s}$ the system of microscopic
equations
\begin{equation}\label{5.10}
\rho_{s}\frac{\partial ^{2}\W}{\partial t^{2}}=
\lambda_{1}\triangle_y \W -\nabla_y R^{s} -\frac{1}{1-m}\nabla_x \pi
+\rho_{s}\F, \quad \y \in Y_{s},
\end{equation}
\begin{equation}\label{5.11}
     \W=\w_{f}, \quad \y \in \gamma
\end{equation}
in the case $\lambda_{1}>0$, and relations
\begin{equation}\label{5.12}
    \rho_{s}\frac{\partial ^{2}\W}{\partial t^{2}}= -\nabla_y R^{s}
    -\frac{1}{1-m}\nabla _{x}\pi
    +\rho_{s}\F, \quad \y \in Y_{s},
\end{equation}
    \begin{equation}\label{5.13}
    (\W-\w_{f})\cdot{\mathbf n}=0, \quad \y \in \gamma
\end{equation}
in the case $\lambda_{1}=0$.

In Eq. (\ref{5.13}) ${\mathbf n}$ is the unit normal to $\gamma$.
\end{lemma}

The proof of this lemma repeats the proof of lemma \ref{lemma4.15}.

In the same way one can prove
\begin{lemma} \label{lemma5.4}
Let $\mu_{1}<\infty$ and $\lambda_{1}=\infty$. Then functions $\{\W,
\w_{s}, \,  q \}$ satisfy in $Y_{f}$ the system of microscopic
equations
\begin{equation}\label{5.14}
\rho_{f}\frac{\partial ^{2}\W}{\partial t^{2}}=
\mu_{1}\triangle_y\frac{\partial\W}{\partial t} -\nabla_y R^{f}
-\frac{1}{m}\nabla_x q +\rho_{f}\F, \quad \y \in Y_{f},
\end{equation}
\begin{equation}\label{5.15}
     \W=\w_{s}, \quad \y \in \gamma
\end{equation}
in the case $\mu_{1}>0$, and relations
\begin{equation}\label{5.16}
    \rho_{f}\frac{\partial ^{2}\W}{\partial t^{2}}= -\nabla_y R^{f}
    -\frac{1}{m}\nabla _{x}q
    +\rho_{f}\F, \quad \y \in Y_{f},
\end{equation}
    \begin{equation}\label{5.17}
    (\W-\w_{s})\cdot{\mathbf n}=0, \quad \y \in \gamma
\end{equation}
in the case $\mu_{1}=0$.
\end{lemma}

\begin{lemma} \label{lemma5.5}
Let $\mu_{1}<\infty$ and $\lambda_{1}<\infty$ and
$\tilde{\rho}=\rho_{f}\chi+\rho_{s}(1-\chi)$. Then functions $\{\W,
 \, \pi\}$ satisfy in $Y$ the system of microscopic
equations
\begin{equation}\label{5.18}
\left. \begin{array}{lll}  \displaystyle &&\tilde{\rho}\partial
^{2}\W /\partial
t^{2}+1 /(1-m)\nabla_{x}\pi-\tilde{\rho}\F=\\[1ex]
&&\div_{y}\{\mu_{1}\chi D(y,\partial\W /\partial t)+
\lambda_{1}(1-\chi)D(y,\W)-R I\}.
 \end{array} \right\}
\end{equation}
\end{lemma}
In the proof of the last lemma we additionally use Nguetseng's
theorem, which states that the sequence $\{\varepsilon
D(x,\partial\w^{\varepsilon}/\partial t)\}$ ($\{\varepsilon
D(x,\w^{\varepsilon})\}$) two-scale converges to
$D(y,\partial\W/\partial t)$ ($D(y,\W)$).\\

\subsection{Homogenized equations}

Lemmas \ref{lemma5.1} and \ref{lemma5.2} imply

\begin{lemma} \label{lemma5.6} Let $\mu_{1}=\lambda_{1}=\infty$, then
$\w_{f}=\w_{s}=\w$ and functions $\w$, $p$, $q$   and $\pi$ satisfy
in $\Omega_{T}$ the system of acoustic equations
\begin{equation} \label{5.19}
\hat{\rho}\frac{\partial ^2\w}{\partial
t^2}=-\frac{1}{(1-m)}\nabla_x\pi+\hat{\rho}\F,
\end{equation}
\begin{equation} \label{5.20}
\frac{1}{p_{*}}p+\frac{1}{\eta_{0}}\pi +\div_x\w=0,
\end{equation}
\begin{equation} \label{5.21}
q=p+\frac{\nu_0}{p_{*}}\frac{\partial p}{\partial t}, \quad
\frac{1}{m}q=\frac{1}{1-m}\pi ,
\end{equation}
homogeneous initial conditions
\begin{equation} \label{5.22}
\w(\x,0)=\frac{\partial \w}{\partial t}(\x,0)=0, \quad \x \in \Omega
\end{equation}
and homogeneous boundary condition
\begin{equation} \label{5.23}
\w(\x,t)\cdot \n(\x)=0, \quad \x\in S,\, t>0.
\end{equation}
\end{lemma}

\begin{lemma} \label{lemma5.7}
Let $\mu_{1}=\infty$ and $\lambda_{1}<\infty$. Then functions $
\w_{f}$, $ \w^{s}$, $p$, $q$  and $\pi$  satisfy in $\Omega_{T}$ the
system of acoustic equations, which consist of the state equations
(\ref{5.21}), balance of momentum equation for the liquid component
\begin{equation} \label{5.24}
\rho_{f}m\frac{\partial ^2\w_{f}}{\partial
t^2}+\rho_{s}\frac{\partial ^2\w^{s}}{\partial
t^2}=-\frac{1}{(1-m)}\nabla_x\pi+\hat{\rho}\F,
\end{equation}
continuity equation
\begin{equation} \label{5.25}
\frac{1}{p_{*}}p+\frac{1}{\eta_{0}}\pi
+m\div_x\w_{f}+\div_x\w^{s}=0,
\end{equation}
 and the relation
\begin{equation}\label{5.26}
\frac{\partial \w^{s}}{\partial t}=(1-m)\frac{\partial
\w_{f}}{\partial t}+\int_{0}^{t} B^{s}_{1}(t-\tau)\cdot
\z^{s}(\x,\tau )d\tau ,
\end{equation}
$$\z^{s}(\x,t)=-\frac{1}{1-m}\nabla_x\pi(\x,t)+\rho_{s}\F(\x,t)-
\rho_{s}\frac{\partial^{2}\w_{f}}{\partial t^{2}}(\x,t)$$ in the
case of $\lambda_{1}>0$, or the balance of momentum equation for the
solid component in the form
\begin{equation}\label{5.27}
\rho_{s}\frac{\partial^{2}\w^{s}}{\partial
t^{2}}=\rho_{s}B^{s}_{2}\cdot \frac{\partial^{2}\w_{f}}{\partial
t^{2}}+((1-m)I-B^{s}_{2})\cdot(-\frac{1}{1-m}\nabla_x\pi+\rho_{s}\F)
\end{equation}
in the case of $\lambda_{1}=0$. The problem (\ref{5.21}),
(\ref{5.24})--(\ref{5.27}) is supplemented by homogeneous initial
conditions (\ref{5.22}) for the displacements in the liquid and the
solid components and homogeneous boundary condition (\ref{5.23}) for
the displacements $\w=m\w_{f}+\w^{s}$.

In Eqs.(\ref{5.26})--(\ref{5.27}) matrices $B^{s}_{1}(t)$ and
$B^{s}_{2}$ are the same as in theorem \ref{theorem2}.
\end{lemma}
\begin{proof}
Eq.(\ref{5.24}) follows directly from Eq.(\ref{5.9}). The continuity
equation (\ref{5.25}) follows  from Eq.(\ref{5.5}) if we take into
account the equality
$$\w=m\w_{f}+\w^{s}.$$
 The
derivation of Eqs. (\ref{5.26})--(\ref{5.27}) is exactly the same as
in lemma \ref{lemma4.7}.
\end{proof}

\begin{lemma} \label{lemma5.8}
Let $\mu_{1}<\infty$ and $\lambda_{1}=\infty$. Then functions $
\w^{f}$, $ \w_{s}$, $p$, $q$ and $\pi$  satisfy in $\Omega_{T}$ the
system of acoustic equations, which consist of the state equations
(\ref{5.21}), the balance of  momentum equation for the solid
component
\begin{equation} \label{5.28}
\rho_{f}\frac{\partial ^2\w^{f}}{\partial
t^2}+\rho_{s}(1-m)\frac{\partial ^2\w_{s}}{\partial
t^2}=-\frac{1}{(1-m)}\nabla_x\pi+\hat{\rho}\F,
\end{equation}
the continuity equation
\begin{equation} \label{5.29}
\frac{1}{p_{*}}p+\frac{1}{\eta_{0}}\pi
+\div_x\w^{f}+(1-m)\div_x\w_{s}=0,
\end{equation}
 and the relation
\begin{equation}\label{5.30}
\frac{\partial \w^{f}}{\partial t}=m\frac{\partial \w_{s}}{\partial
t}+\int_{0}^{t} B^{f}_{1}(t-\tau)\cdot \z^{f}(\x,\tau )d\tau ,
\end{equation}
$$\z^{f}(\x,t)=-\frac{1}{m}\nabla_x q(\x,t)+\rho_{f}\F(\x,t)-
\rho_{f}\frac{\partial^{2}\w_{s}}{\partial t^{2}}(\x,t)$$ in the
case of $\lambda_{1}>0$, or the balance of momentum equation for the
liquid component in the form
\begin{equation}\label{5.31}
\rho_{f}\frac{\partial^{2}\w^{f}}{\partial
t^{2}}=\rho_{f}B^{f}_{2}\cdot \frac{\partial^{2}\w_{s}}{\partial
t^{2}}+(mI-B^{f}_{2})\cdot(-\frac{1}{m}\nabla_x q+\rho_{f}\F)
\end{equation}
in the case of $\lambda_{1}=0$. The problem (\ref{5.21}),
(\ref{5.28})--(\ref{5.31}) is supplemented by homogeneous initial
conditions (\ref{5.22}) for the displacements in the liquid and the
solid components and homogeneous boundary condition (\ref{5.23}) for
the displacements $\w=\w^{f}+(1-m)\w_{s}$.

In Eqs.(\ref{5.30})--(\ref{5.31}) matrices $B^{f}_{1}(t)$ and
$B^{f}_{2}$ are given below by formulas (\ref{5.32})--(\ref{5.33}).
\end{lemma}
\begin{proof}
The proof of this lemma repeats  proofs of previous lemmas and
\begin{equation}\label{5.32}
B^{f}_{1}(t)= \langle \sum_{i=1}^{3}\V^{i}(\y,t)\rangle
_{Y_{f}}\otimes {\mathbf e}_{i},
\end{equation}
\begin{equation}\label{5.33}
B^{f}_{2}=\sum_{i=1}^{3}\langle \nabla R^{f}_{i}(\y)\rangle
_{Y_{f}}\otimes {\mathbf e}_{i},
\end{equation}
 where functions $\V^{i}(\y,t)$ solve the periodic
initial-boundary value problem
\begin{equation}\label{5.34}
\left. \begin{array}{lll}  \displaystyle \rho_{f}\frac{\partial
\V^{i}}{\partial t}-\mu_{1}\triangle \V^{i} +\nabla R^{i} =0,
  \quad \mbox{div}_y \V^{i} =0, \quad \y \in Y_{f},  t>0,\\[1ex]
\V^{i}=0, \quad \y \in \gamma ,  t>0;\quad
 \rho_{f}\V^{i}(y,0)={\mathbf e}_{i}, \quad \y \in Y_{f},
 \end{array} \right\}
\end{equation}
and functions $R^{f}_{i}(\y,t)$ solve the periodic boundary value
problem
\begin{equation}\label{5.35}
\triangle_y R^{f}_{i}=0,\quad \y \in Y_{f}; \quad \nabla_y
R^{f}_{i}\cdot \n =\n\cdot {\mathbf e}_{i}, \quad \y \in \gamma .
\end{equation}
Note that as before, the matrix $((mI-B^{f}_{2})$ is symmetric and
positively defined.
\end{proof}
The proof of theorem \ref{theorem3} is completed by
\begin{lemma} \label{lemma5.9}
Let Let $\mu_{1}<\infty$ and $\lambda_{1}<\infty$. Then functions $
\w$, $p$, $q$  and $\pi$  satisfy in $\Omega_{T}$ the system of
acoustic equations, which consist of the continuity and the state
equations (\ref{5.20})  and (\ref{5.21}) and the relation
\begin{equation}\label{5.36}
\frac{\partial\w}{\partial
t}=\int_{0}^{t}B^{\pi}(t-\tau)\cdot\nabla\pi (\x,\tau )d\tau
+\textbf{f}(\x,t),
\end{equation}
where  $B^{\pi}(t)$ and $\textbf{f}(\x,t)$ are given below by
Eqs.(\ref{5.40}) and (\ref{5.41}).

The problem (\ref{5.20}), (\ref{5.21}), (\ref{5.36}) is supplemented
by homogeneous initial and boundary conditions (\ref{5.22}) and
(\ref{5.23}).
\end{lemma}
\begin{proof}
Let
$$\W=\int_{0}^{t}\sum_{i=1}^{3}\{\W^{\pi}_{i}(\y,t-\tau)\frac{\partial\pi}{\partial x_{i}}(\x,\tau)+
\W^{F}_{i}(\y,t-\tau)F_{i}(\x,\tau)\}d\tau,$$
$$R=\int_{0}^{t}\sum_{i=1}^{3}\{R^{\pi}_{i}(\y,t-\tau)\frac{\partial\pi}{\partial x_{i}}(\x,\tau)+
R^{F}_{i}(\y,t-\tau)F_{i}(\x,\tau)\}d\tau ,$$ where
$\F=\sum_{i=1}^{3}F_{i}{\mathbf e}_{i}$  and functions
$\{\W^{\pi}_{i}(\y,t),R^{\pi}_{i}(\y,t)\}$ and
$\{\W^{F}_{i}(\y,t),R^{F}_{i}(\y,t)\}$ are periodic in $\y$
solutions  of the system
\begin{equation}\label{5.37}
\left. \begin{array}{lll}  \displaystyle \div_{y}\{\mu_{1}\chi
D(y,\partial\W^{j}_{i}
/\partial t)+ \lambda_{1}(1-\chi)D(y,\W^{j}_{i})-R^{j}_{i}I\}=\\[1ex]
\tilde{\rho}\partial ^{2}\W^{j}_{i}/\partial t^{2},\quad
div_{y}\W^{j}_{i}=0, \quad \y \in Y, \quad t>0, \quad j=\pi, F,
 \end{array} \right\}
\end{equation}
which satisfy the following initial conditions
\begin{equation} \label{5.38}
\W^{\pi}_{i}(\y,0)=0, \quad \tilde{\rho}\frac{\partial
\W^{\pi}_{i}}{\partial t}(\y,0)=-\frac{1}{1-m}{\mathbf e}_{i}, \quad
\x \in Y,
\end{equation}
\begin{equation} \label{5.39}
\W^{F}_{i}(\y,0)=0, \quad \frac{\partial \W^{F}_{i}}{\partial
t}(\y,0)={\mathbf e}_{i}, \quad \x \in Y.
\end{equation}
Then the functions $\W$ and $R$ solve the system of microscopic
equations (\ref{5.7}) and (\ref{5.18})  and by definition
$\w=\langle\W\rangle_{Y}$. Therefore
\begin{equation}\label{5.40}
B^{\pi}(t)= \sum_{i=1}^{3}\langle \frac{\partial
\W^{\pi}_{i}}{\partial t}(\y,t)\rangle _{Y}\otimes{\mathbf e}_{i},
\end{equation}
\begin{equation}\label{5.41}
\textbf{f}(\x,t)= \int_{0}^{t}\sum_{i=1}^{3}\langle \frac{\partial
\W^{F}_{i}}{\partial t}\rangle _{Y}(t-\tau)F_{i}(\x,\tau)d\tau .
\end{equation}
The solvability and the uniqueness of  problems (\ref{5.37}),
(\ref{5.38}) or (\ref{5.37}), (\ref{5.39}) follow directly from the
energy identity
$$\frac{1}{2}\langle \tilde{\rho}(\frac{\partial \W^{j}_{i}}{\partial
t})^{2}\rangle _{Y}(t)+\frac{1}{2}\langle
\lambda_{1}D(y,\W^{j}_{i}):D(y,\W^{j}_{i})\rangle _{Y_{s}}(t)+$$
$$ \int_{0}^{t}\langle \mu_{1}D(y,\frac{\partial\W^{j}_{i}}{\partial \tau}):
D(y,\frac{\partial\W^{j}_{i}}{\partial \tau})\rangle
_{Y_{s}}(\tau)d\tau=\frac{1}{2}\beta^{j}, $$ for $i=1,2,3$ and
$j=\pi, F$.

Here
$$\beta^{\pi}=\langle \frac{1}{\tilde{\rho}}\rangle _{Y},
 \quad \beta^{F}=\langle \tilde{\rho}\rangle _{Y}.$$
 As before, equations (\ref{5.37}) are understood in the sense of
 distributions and the function $B^{\pi}(t)$ has no time derivative at
 $t=0$. That is why  we cannot represent relation (\ref{5.36}) in
 the form of the balance of momentum equation, like (\ref{5.19}) or
 (\ref{5.27}).
\end{proof}

\end{document}